\documentclass[preprint,review,3p,12pt]{elsarticle}

\usepackage{amssymb,dsfont}
\usepackage{lineno}
\usepackage{graphicx}
\usepackage{epstopdf}
\usepackage{underscore}
\usepackage{bm}
\usepackage{color}
\usepackage[dvipsnames]{xcolor}

\newtheorem{definition}{Definition}
\newtheorem{principle}{Principle}
\biboptions{comma,sort}
\journal{Journal of Computational Physics}
\date{}

\usepackage{soul}

\begin{document}

\begin{frontmatter}
  \title{An RBF-FD closest point method for solving PDEs on surfaces}
  \author[rvt,rvt2]{A.~Petras\corref{cor1}}
  \ead{apetras@bcamath.org}
  \author[rvt3]{L.~Ling}
  \author[rvt]{S.J.~Ruuth}
  \ead{sruuth@sfu.ca}
  \cortext[cor1]{Corresponding author}
  \address[rvt]{Department of Mathematics, Simon Fraser University, Burnaby, British Columbia, Canada V5A1S6}
  \address[rvt2]{BCAM-Basque Center for Applied Mathematics, Bilbao, Basque Country, Spain 48015}
  \address[rvt3]{Department of Mathematics, Hong Kong Baptist University, Kowloon Tong, Hong Kong}

  \begin{abstract}
    Partial differential equations (PDEs) on surfaces appear in many applications throughout the natural and applied sciences.
    The classical closest point method (Ruuth and Merriman, J. Comput. Phys. 227(3):1943-1961, [2008]) is an embedding method for solving PDEs on surfaces using standard finite difference schemes.
In this paper, we formulate an explicit closest point method using finite difference schemes derived from radial basis functions (RBF-FD). Unlike the orthogonal gradients method (Piret, J. Comput. Phys. 231(14):4662-4675, [2012]), our proposed method uses RBF centers on regular grid nodes. This formulation not only reduces the computational cost but also avoids the ill-conditioning from point clustering on the surface and is more natural to couple with a grid based manifold evolution algorithm (Leung and Zhao, J. Comput. Phys. 228(8):2993-3024, [2009]).
When compared to the standard finite difference discretization of the closest point method, the proposed method requires a smaller computational domain surrounding the surface, resulting in a decrease in the number of sampling points on the surface. In addition, higher-order schemes can easily be constructed by increasing the number of points in the RBF-FD stencil. Applications to a variety of examples are provided to illustrate the numerical convergence of the method.
\end{abstract}

  \begin{keyword}
  closest point method \sep embedding method \sep radial basis functions \sep finite differences
  \end{keyword}
\end{frontmatter}

\section{Introduction}
Many applications in the natural and applied sciences require the solution of partial differential equations (PDEs) on surfaces. Image processing applications include the placement of an image on a surface \cite{turk1991generating},
the restoration of a damaged pattern on a surface \cite{bertalmio2001navier} and the segmentation and the denoising of images on surfaces \cite{tian2009segmentation,biddle2013volume}. In biology, applications include the formation of patterns on animal coats \cite{murray2001mathematical} and the wound healing process \cite{olsen1998spatially}. In computer graphics, applications are found
in the topic of real time fluid visualization on surfaces \cite{auer2012real}.

Various numerical methods have been developed to approximate the solution of PDEs on surfaces. These include methods applied on parametrized surfaces, on triangulated surfaces and on surfaces embedded in a higher dimensional space. Solution of PDEs on parametrized surfaces can be efficient for surfaces where a parametrization is possible \cite{lui2005solving,floater2005surface}, however a parametrization of a
surface often leads to distortions of the surface and singularities \cite{floater2005surface}. Triangulated surfaces
avoid these issues. Finite difference methods can be applied to solve PDEs on triangulated surfaces \cite{turk1991generating}, but there are difficulties in the calculation of geometric quantities, including the normal vector and the curvature of a surface \cite{bertalmio2001variational}. On the other hand, finite element methods on triangulated surfaces are effective in solving parabolic or elliptic PDEs \cite{dziuk2007surface}. Methods using surfaces embedded into a $d$-dimensional space extend the surface PDE in the embedding space and solve the extended PDE using standard Cartesian methods. 

A popular method employs a level set representation of surfaces and  a projection operator to solve surface PDEs \cite{greer2006improvement,flyer2009radial}.
Typically, the computational domain consists of points in a neighborhood of the surface. This may lead to the introduction of artificial boundary conditions
at the boundary of the computational domain which can degrade the accuracy of the method
\cite{cheung2015localized}. Meshfree approximations using radial basis functions (RBFs) are also becoming popular within the embedded surfaces class \cite{fuselier2013high,shankar2015radial}.

The closest point method \cite{ruuth2008simple} is an embedding method which uses a closest point representation of the surface to solve PDEs on surfaces. In the classical formulation \cite{ruuth2008simple,macdonald2008level},
the discretization is carried out in a neighborhood of the surface using standard finite difference schemes and barycentric Lagrangian interpolation.
The implicit closest point method was introduced in \cite{macdonald2009implicit} to provide a stable approximation of the Laplace-Beltrami and other higher-order surface differential operators. Application of the implicit closest point method to the solution of eigenvalue problems appears in \cite{macdonald2011solving}. See also \cite{marz2012calculus} for a study of the theoretical foundation of the closest point method.

The extension via the closest point mapping is also used as part of the development for other methods for solving surface PDEs. In \cite{piret2012orthogonal}, the author uses the closest point mapping to derive an RBF method for solving surface PDEs.
The method gives a high-order approximation to the solutions of surface PDEs in a variety of examples. See also  \cite{cheung2015localized} for a related RBF method that carries out
a local approximation of surface differential operators to solve PDEs on folded surfaces. In addition, the extension via the closest point mapping is used in the computation of integrals over curves and surfaces \cite{kublik2016integration} and in the solution of PDEs on closed, smooth surfaces using volumetric variational principles \cite{chu2017volumetric}.

In this paper, we propose a new method using a closest point representation of a surface and
finite difference stencils derived from radial basis functions (RBF-FD).
Notably, the use of RBF-FD leads to a method (RBF-CPM) that evaluates derivatives on the surface, rather than in the embedding space. This eliminates the interpolation step in the
evaluation of derivatives, thereby eliminating a potential source of error and computational cost. In addition, standard RBF-FD and global RBF methods for surface PDEs may suffer ill-conditioning due to small separating distances in the surface points \cite{wendland2004scattered}. Our method uses RBF-FD stencils on regular Cartesian grid nodes, thus allowing irregular collocation points on the surface and avoiding the ill-conditioning that may arise due to point clustering on the surface. Due to the regularity of the RBF-FD stencil, the collocation matrix associated with the calculation of the RBF-FD weights is independent of the surface, and its inverse can be accurately calculated locally. Due to repeated patterns on the RBF centers, only a small number of collocation matrix inverses need to be calculated, thus reducing the computational cost over existing RBF methods.

In our method, second-order accuracy in $\Delta x$ can be achieved with smaller computational domains and fewer points on the surface
relative to the classical closest point method. Furthermore, higher-order schemes are obtained simply by increasing the number of points in the finite difference stencil.

The paper unfolds as follows. In Section~\ref{MethodsReviewSection}, we review the classical closest point method and RBF approximation.
Section~\ref{CPMRBFFD} gives our new method
and studies the selection of parameters and the computational domain.
Section~\ref{Numerics} considers the performance of the method using a variety of convergence studies and numerical examples
in two and three dimensions.  Finally Section~\ref{SummarySection}  concludes and explores potential future work.

\section{Numerical methods review}\label{MethodsReviewSection}

\subsection{The closest point method}\label{CPM}
The classical closest point method \cite{ruuth2008simple} is a
simple numerical method for approximating the solution of PDEs on surfaces.
In this section, we review the method and its components.
\subsubsection{Surface representation}\label{section: CPCalculation}
To begin, the {\it closest point} to the surface is defined:
\begin{definition}\label{ClosestPointRepresentation}
  Let $\Gamma\subset\mathbb{R}^d$ be a surface and $\mathbf{z}\in \Omega\subset \mathbb{R}^d$ be some point in the embedding space $\Omega\supset\Gamma$. 
  Then,
  $$cp_\Gamma(\mathbf{z})=\arg\min_{\mathbf{x}\in \Gamma}\|\mathbf{x}-\mathbf{z}\|_2$$
  is the closest point of $\mathbf{z}$ to the surface $\Gamma$.
\end{definition}
In a neighborhood of the surface,  $cp_\Gamma$ will be $C^p$-smooth for
a $C^{p+1}$-smooth surface $\Gamma$ \cite{marz2012calculus}.

To discretize, a Cartesian grid is introduced in the embedding space,
over a neighborhood of the surface.
Typically, this neighborhood includes all grid nodes whose Euclidean distance to the
surface is less than or equal to some constant $\gamma_{CPM}$.
Following a common convention (e.g., \cite{leung2009grid}),
we refer to this localized computational domain as the {\it computational tube}, and
the corresponding radius  $\gamma_{CPM}$ as the {\it computational tube radius}.  
To avoid introducing discontinuities into $cp_\Gamma$,
the computational tube radius should satisfy $\gamma_{CPM}<\kappa_\infty^{-1}$, where $\kappa_\infty$ is an upper
bound on the curvatures of $\Gamma$ \cite{chu2017volumetric}.
The grid points and their closest points together form a {\it closest point representation} of the surface $\Gamma$.

The method used to determine the closest point function depends on the type of surface under consideration.
For simple surfaces, such as the sphere and the torus, an analytical formula for the closest point function is available, and
the preferred approach is to simply evaluate the formula.   On the other hand, for ellipsoids, the M\"{o}bius strip,
and many other interesting shapes, the surface may be given in parameterized form.
Here, standard numerical optimization techniques can be applied to find the closest point on the surface (cf. \cite{merrimanruuth2007}).
Finally, we consider surfaces in triangulated form.  In this case, we follow  \cite{macdonald2008level}
and loop over the list of triangles.
For grid nodes near a triangle $T_i$ (specifically, nodes that are within a Euclidean distance $\gamma_{CPM}$ of $T_i$),
we compute and store the closest point on $T_i$.  For each grid node, the closest point over all stored possibilities
is the closest point on the surface.  See \cite{macdonald2008level} for further details on this procedure. Another approach which follows the causality of the eikonal equation can
also be applied \cite{tsai2002rapid}.

\subsubsection{The equivalence principles}
In the closest point method, we do not solve the surface PDE problem on the surface directly.
Instead, we solve a suitable differential equation defined over the computational domain.
The key property of this differential equation (the {\it embedding equation}) is that its solution
{\it on the surface} must agree with the solution to the original surface PDE.
Values off the surface do not directly give the solution to the PDE-on-surface problem.

To form the embedding equation, surface derivatives are replaced with closest point operators
and standard Cartesian derivatives.
  The foundation of this procedure consists of two principles \cite{ruuth2008simple,cheung2017Kernel}:
the \emph{equivalence of gradients} and the \emph{equivalence of divergence}.
\begin{principle}\label{equivalence of gradients}
  Let $v$ be any function on $\mathbb{R}^d$ that is constant along normal directions of $\Gamma$. Then, at the surface, intrinsic gradients are equivalent to standard gradients, $\nabla_\Gamma v=\nabla v$.
\end{principle}

\begin{principle}\label{equivalence of divergence}
  Let $\mathbf{v}$ be any vector field on $\mathbb{R}^d$ that is tangent to $\Gamma$ and tangent to all surfaces displaced by a fixed distance from $\Gamma$. Then, at the surface, $\nabla_\Gamma\cdot\mathbf{v}=\nabla\cdot\mathbf{v}$.
\end{principle}

General surface differential operators can be replaced with the corresponding Cartesian differential operators by combining the two principles.
Of particular interest to us is the composition of the divergence and gradient operators:
Let $u$ be any function on $\mathbb{R}^d$ that is constant along normal directions of $\Gamma$. Then, on the surface, the Laplace-Beltrami operator is equivalent to the standard Laplacian operator,
$\Delta_\Gamma u=\Delta u$ for all $\mathbf{x}\in\Gamma$.
This property is referred to as the \emph{equivalence of the Laplacian}.
As a consequence, heat flow intrinsic to a surface can be approximated by alternating constant-along-normal extension
with standard heat flow in the underlying embedding space.
Other, similarly straightforward combinations of the gradient and divergence principles lead to replacements for
nonlinear diffusive flows such as curvature motion intrinsic to a surface \cite{merrimanruuth2007,ruuth2008simple}.
Higher-order operators can also be approximated in the embedding space, although additional extension operators may be needed.
See, for example, \cite{macdonald2009implicit} where two extension operators are used as part of the
approximation of fourth-order operators.

\subsubsection{The closest point method}
Evolution of the embedding equation may be carried out in a similar fashion, to yield the explicit closest point method \cite{ruuth2008simple}.
Specifically, given a closest point representation of a surface $\Gamma$, the explicit closest point method alternates between the following two steps:
\begin{enumerate}
  \item \textbf{Closest point extension}.  Carry out a constant-along-normal extension
  of $u:\Gamma\to\mathbb{R}$ to yield $\tilde{u}:\Omega\to\mathbb{R}$ by $\tilde{u}(\mathbf{z})=u(cp_\Gamma(\mathbf{z}))$
  for each $\mathbf{z}$ in the tubular computational domain $\Omega\supset\Gamma$.
  \item \textbf{Evolution}.  The PDE is solved on the tubular computational domain $\Omega$ in the embedding space for one time step (or one stage of a Runge-Kutta method).
\end{enumerate}

Note that the closest point extension step is an interpolation step since $cp_\Gamma(\mathbf{z})$ is not necessarily a grid point.
In \cite{ruuth2008simple}, barycentric Lagrange interpolation is used with polynomial degree $p=q+r-1$, where $q$ is the order of finite differences schemes and $r$ is the order of the derivatives.
Localization of the computation is accomplished by computing over a computational tube surrounding the
surface \cite{ruuth2008simple,macdonald2009implicit}.
For second-order finite differences and second-order derivatives, it is sufficient to choose a
computational tube radius of $\gamma_{CPM}$ where
\begin{equation}\label{Bandwidth}
  \gamma_{CPM}=\sqrt{(d-1)\left(\frac{p+1}{2}\right)^2+\left(1+\frac{p+1}{2}\right)^2}\Delta x
\end{equation}
in a $d$-dimensional embedding space \cite{ruuth2008simple}.

\subsection{Global Radial Basis Function (RBF) approximation}\label{sec2.2}
RBF approximation is a powerful tool for approximating smooth functions on a variety of geometries. Following \cite{fasshauer2007meshfree,fornberg2015solving},
given an RBF $\phi(r)$ (see Table~\ref{RBFs} for some RBF options) and a set of scattered
points $Z = \{\mathbf{z}_j\}_{j=1}^{n_Z}$ called RBF centers, the RBF interpolant has the form
\begin{equation}\label{RBFInterpolant}
s(\mathbf{x})= \sum_{j=1}^{n_Z}\lambda_j\phi(\|\mathbf{x}-\mathbf{z}_j\|)
\end{equation}
with coefficients $\lambda_j$. For the RBF interpolation of any smooth function $v:\Omega\to\mathbb{R}$, the coefficients $\lambda_j$ can be found by interpolation conditions at $Z\subset\Omega$, i.e. by solving the symmetric linear system
\begin{equation}\label{RBFInterpolantMatrixForm}
\displaystyle A(Z,Z) \bm{\lambda} = v(Z),
\end{equation}
for $\bm{\lambda}=[\lambda_j]\in \mathbb{R}^{n_Z \times 1}$, where ${v}(Z):=[v(\mathbf{z}_j)]\in \mathbb{R}^{n_Z \times 1}$ and $A(Z,Z):=[\phi(\|\mathbf{z}_i-\mathbf{z}_j\|)]\in \mathbb{R}^{n_Z\times n_Z}$ for $\mathbf{z}_i,\mathbf{z}_j\in Z$ in an orderly sense.

\begin{table}[h]
\centering
\begin{tabular}{lll}
  \hline
  Name of RBF & Abbreviation & Definition \\ \hline
  \underline{Smooth RBFs} & & \\
  Gaussian & GA & $\phi(r) = e^{-(\epsilon r)^2}$\\
  Multiquadratic & MQ & $\phi(r)=\sqrt{1+(\epsilon r)^2}$ \\
  Inverse multiquadratic & IMQ &  $\displaystyle\phi(r)=\frac{1}{\sqrt{1+(\epsilon r)^2}}$\\
  Inverse quadratic & IQ & $\displaystyle\phi(r)= \frac{1}{1+(\epsilon r)^2}$\\
&&\\
  \underline{Piecewise smooth RBFs} & & \\
  Cubic & CU & $\phi(r)=|r|^3$ \\
  Thin plate spline & TPS & $\phi(r)=r^2\ln|r|$ \\
  \hline
\end{tabular}
\caption{Definition of some commonly used RBFs.}\label{RBFs}
\end{table}
Then, we can use the interpolant to approximate the derivatives of $v$.
If $L$ is a differential operator, then the quantity $Lv$ at some point $\mathbf{x}$ can be approximated as
\begin{equation}\label{RBFDer}
  Lv(\mathbf{x})\approx Ls(\mathbf{x}) = \sum_{j=1}^{n_Z}\lambda_jL\phi(\|\mathbf{x}-\mathbf{z}_j\|),
\end{equation}
in which $L$ acts upon the variable $\mathbf{x}$ in the basis function $\phi$.
In matrix form, we can express (\ref{RBFDer}) as
\begin{equation}\label{RBFDerU}
  Lv(\mathbf{x})\approx B(\mathbf{x},Z)A(Z,Z)^{-1} v(Z),
\end{equation}
where $B(\mathbf{x},Z) = [L\phi(\|\mathbf{x}-\mathbf{z}_1\|),\ldots,L\phi(\|\mathbf{x}-\mathbf{z}_{n_Z}\|)]\in \mathbb{R}^{1\times n_Z}$.
This expression provides a radial basis function
pseudo-spectral method, which can be easily localized to give a RBF
finite difference discretization (RBF-FD) for the differential operator $L$ evaluated at the data site $\mathbf{x}\in\Omega$.
Some compact RBF-FD stencils can be found in \cite{fornberg2015primer} for applications in geosciences. Details on RBF-FD stencils that use a given number of nearest neighbors are found in \cite{flyer2016role}.

\section{A closest point method for solving PDEs on surfaces using RBF-FD}\label{CPMRBFFD}
In this section, we introduce an explicit RBF closest point method for solving PDEs on surfaces.
Our method uses RBF-FD stencils that consist of $m$ closest neighboring grid points.
As part of our method, we provide an approach for the calculation of a computational tube around the surface.

For illustration purposes, consider the heat equation
\begin{equation}\label{equation:surface heat}
\displaystyle u_t = \Delta_\Gamma u
\end{equation}
intrinsic to a surface $\Gamma$, where $u:\Gamma\to\mathbb{R}$ is a function defined solely on the surface. Using a set of specially constructed surface data points $X = \{\mathbf{x}_j\}\subset \Gamma$, our aim is to design an efficient, accurate, and robust finite difference scheme in order to spatially discretize equation~(\ref{equation:surface heat}), i.e.
$$[\Delta_\Gamma u](X)\approx W u(X),$$
where $W$ is a
differentiation matrix that takes the vector of function values $u(X) := [u(\mathbf{x}_j)]_{\mathbf{x}_j\in X}$ to the approximated values of  $\Delta_\Gamma u$ at $X$.

In the literature, this can be done by the orthogonal gradients method \cite{piret2012orthogonal} and projection methods \cite{fuselier2013high} using RBFs, in which the geometry of $X$ plays a key role in the numerical stability. In the orthogonal gradients method, the RBF centers $Z = \{\mathbf{z}_j\}$ are constructed by extending the surface points $X$ in the embedding space in the normal direction. The geometry of the surface $\Gamma$ and the points $X$ determine that of $Z$, which should ideally be quasi-uniform. Projection methods work solely on $X$
with the corresponding analysis carried out in Sobolev spaces on $\Gamma$. Like other typical kernel approximation theories, convergence comes when $X$ gets dense, i.e. as the fill distance $h_X\longrightarrow0$. Both approaches involve solving interpolation problems, whose conditioning depends on the minimum separating distance $q_X$ of $X$ \cite{wendland2004scattered}. Thus, it is common to require that the mesh ratio $\rho_X:=h_X/q_X\geq1$ of the surface points is bounded. In short, quasi-uniform points $X$ need to be used on surfaces, which may not be an easy task.

Finite difference schemes, on the other hand, work on regular grids with mesh ratio exactly equal to 1. Yet, extra work is required to apply finite differences to surfaces in general (e.g., the original closest point method \cite{ruuth2008simple}). We propose a new RBF kernel based formulation to get the best of both worlds.

\subsection{Description of the method}\label{Method description}
\begin{figure}[h!]
\centering
    \includegraphics[width=0.49\textwidth]{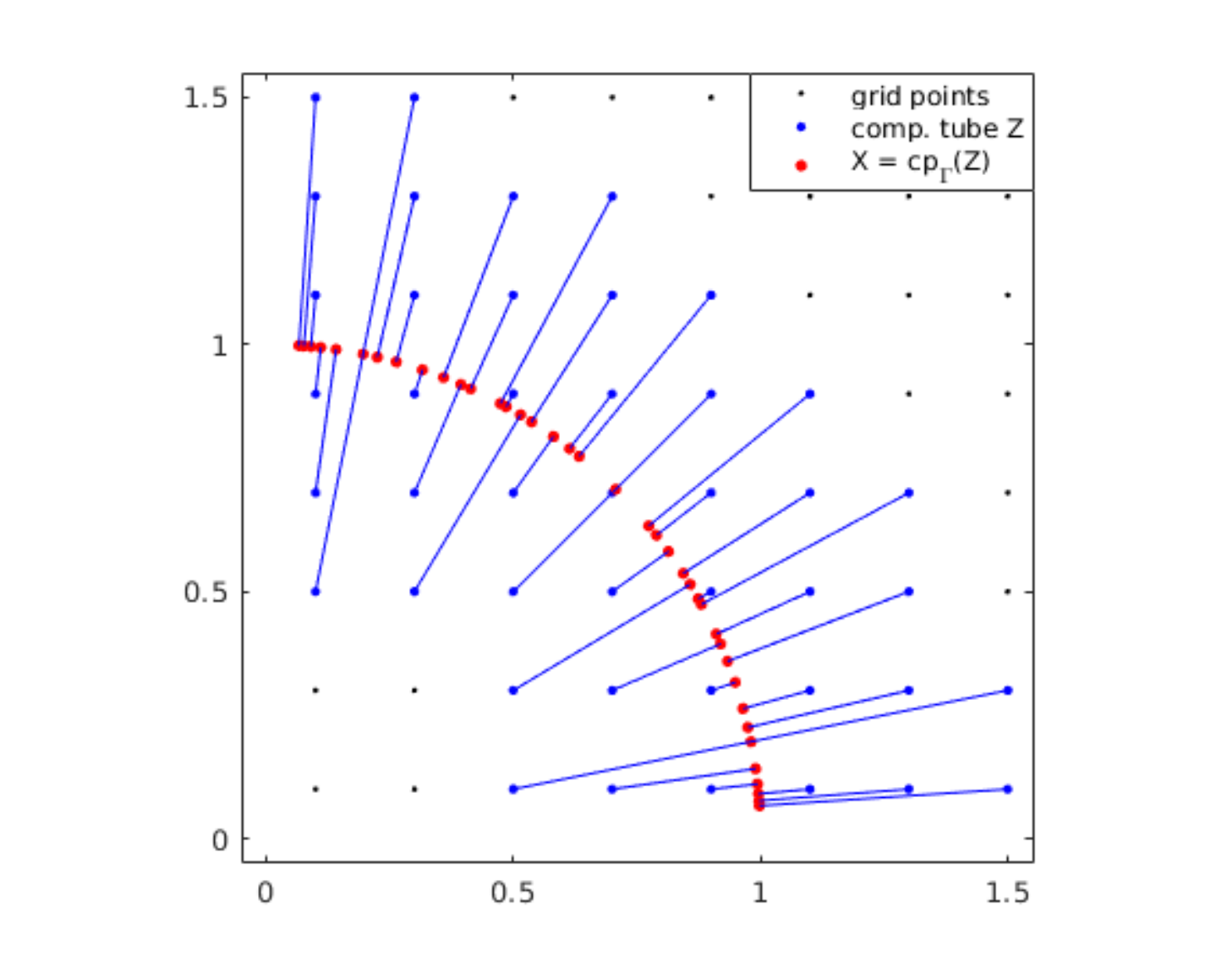}
    \includegraphics[width=0.49\textwidth]{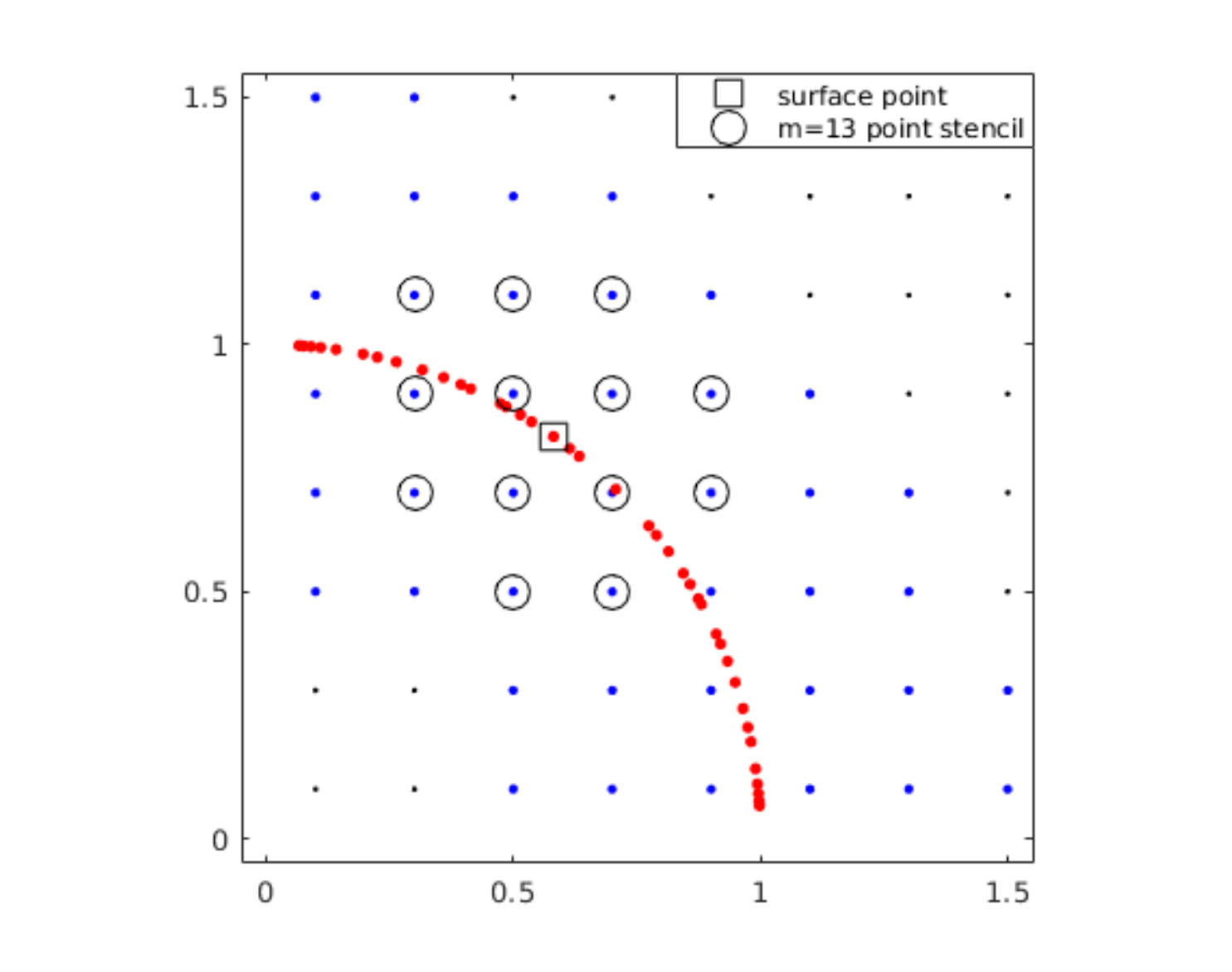}
    \caption{Left: The closest points (red dots) to the grid points in the computational domain (blue dots) on the surface. Right: An example of a $m=13$ point stencil (circled blue dots) for a surface point (squared red dot).}
    \label{RBFStencil}
\end{figure}

We start with a collection of  Cartesian grid points $Z = \{\mathbf{z}_j\}_{j=1}^{n_Z}\subset\Omega$ on a small tubular domain containing the surface $\Gamma$. Then, we define surface data points via $\mathbf{x}_j = cp_\Gamma(\mathbf{z}_j)$ to form a set $X=\{\mathbf{x}_j\}_{j=1}^{n_Z}\subset \Gamma$;
see Figure~\ref{RBFStencil} for a schematic demonstration.
Applying the equivalence of the Laplacian property yields the relation
$$
\Delta_\Gamma u(X) = \Delta \tilde{u}(X),
$$
which holds for any $X\subset\Gamma$, and where we denote the constant-along-normal extension of  $u$ by $\tilde{u}:\Omega\to\mathbb{R}$.

We now deploy the methodology locally to obtain an RBF-FD approximation.
For each surface point $\mathbf{x}_j = cp_\Gamma(\mathbf{z}_j)$ for some $\mathbf{z}_j\in Z$ and $j=1,\ldots,n_Z$, let  $Z_j=\{\mathbf{z}_{j_1},\ldots,\mathbf{z}_{j_m}\}\subset Z$ denote the $m$ nearest neighborhood of $\mathbf{x}_j$. Locally, we take the $\mathbf{x}_j$--local interpolant of $\tilde{u}$ using basis function $\phi$ at centers $Z_j$, denoted by $s_j$ below, as the function $v$ in Section \ref{sec2.2}.
Then, (\ref{RBFDerU}) gives an approximation scheme
\begin{eqnarray*}
  \Delta \tilde{u}(\mathbf{x}_j) &\approx&  \Delta  s_j (\mathbf{x}_j)\\
  &=& B(\mathbf{x}_j,Z_j)A(Z_j,Z_j)^{-1} \tilde{u}(Z_j) \\
  &=:& \mathbf{w}_j \tilde{u}(Z_j).
\end{eqnarray*}
In other words, the nonzero RBF-FD weight $\mathbf{w}_j\in\mathbb{R}^{1\times m}$ is given as the product of a row vector
$B(\mathbf{x}_j,Z_j) = [\Delta\phi(\|\mathbf{x}_j-\mathbf{z}_{j_1}\|),\ldots,\Delta\phi(\|\mathbf{x}_j-\mathbf{z}_{j_m}\|)]\in \mathbb{R}^{1\times m}$ and the inverse matrix of $A(Z_j,Z_j)=[\phi(\|\mathbf{z}_{j_k}-\mathbf{z}_{j_\ell}\|)]\in\mathbb{R}^{m\times m}$ for $\mathbf{z}_{j_k},\mathbf{z}_{j_\ell}\in Z_j$.
Using $\mathbf{w}_j$ for $j=1,\ldots,n_Z$, we can assemble the RBF-FD matrix $W$ such that
\begin{equation}\label{Wu}
  \Delta_\Gamma u(X) =  \Delta \tilde{u}(X) \approx W \tilde{u}(Z).
\end{equation}

In addition, we also construct an RBF-FD projection matrix $P$ such that
\begin{equation}\label{Pu}
 u(X) = \tilde{u}(X) \approx P \tilde{u}(Z),
\end{equation}
simply by replacing $\Delta$ with the identity map in the computation.   The matrix $P$ will appear later as part of our time-stepping scheme.

Note that all approximations are done in the embedding space, which is independent of the geometry of $X$. By using a regular $Z$, we obtain
collocation matrices $A(Z_j,Z_j)$ that depend on the geometry of $Z_j$, but are independent of the surface $\Gamma$. Thus, we can precompute all $A(Z_j,Z_j)^{-1}$, and the evaluation of each RBF-FD weight requires only matrix-vector multiplications. As a consequence,
we can employ accurate and expensive solvers to compute the $A(Z_j,Z_j)^{-1}$ matrices without harming the overall performance of the proposed method. This greatly improves computational efficiency over the existing methods in which no matrix inverse can be reused due to the lack of repeated pattern in data point geometry.

The RBF finite difference stencil used in this paper consists of the $m$ closest grid points $Z_j$, to each surface point $\mathbf{x}_j$. Methods that use RBF-FD stencils of a given number of nearest neighbors on scattered nodes require the use of quasi-uniform nodes on the surface \cite{shankar2015radial}. By using RBF-FD stencils on the grid nodes, we avoid the ill-conditioning of the RBF-FD matrices that arises from the clustering of nodes on the surface. Figure~\ref{RBFStencil} shows an example of an RBF-FD stencil that consists of the $m=13$ closest grid points to a particular surface point (displayed using a squared red dot).

We are now ready to temporally discretize the heat equation
\[
    {u}_t = \Delta_\Gamma {u}
\]
intrinsic to the surface $\Gamma$.
Using the forward Euler scheme with spatial  discretization at $X\subset \Gamma$ as above, we have
\begin{equation}\label{forwardEuler}
     {u}(X,t^{n+1}) = {u}(X,t^{n})  + \Delta t \Delta_\Gamma {u}(X,t^n) + \mathcal{O}(\Delta t^2),
\end{equation}
for $t^n=n\,\Delta t$.
Recall that in our proposed setup $Z\in\Omega$ is regular, whereas $X=cp_\Gamma(Z)$ could be highly nonuniform.
It is desirable to work on the $Z$ nodes.

\begin{figure}[h!]
\centering
    \includegraphics[width=0.49\textwidth]{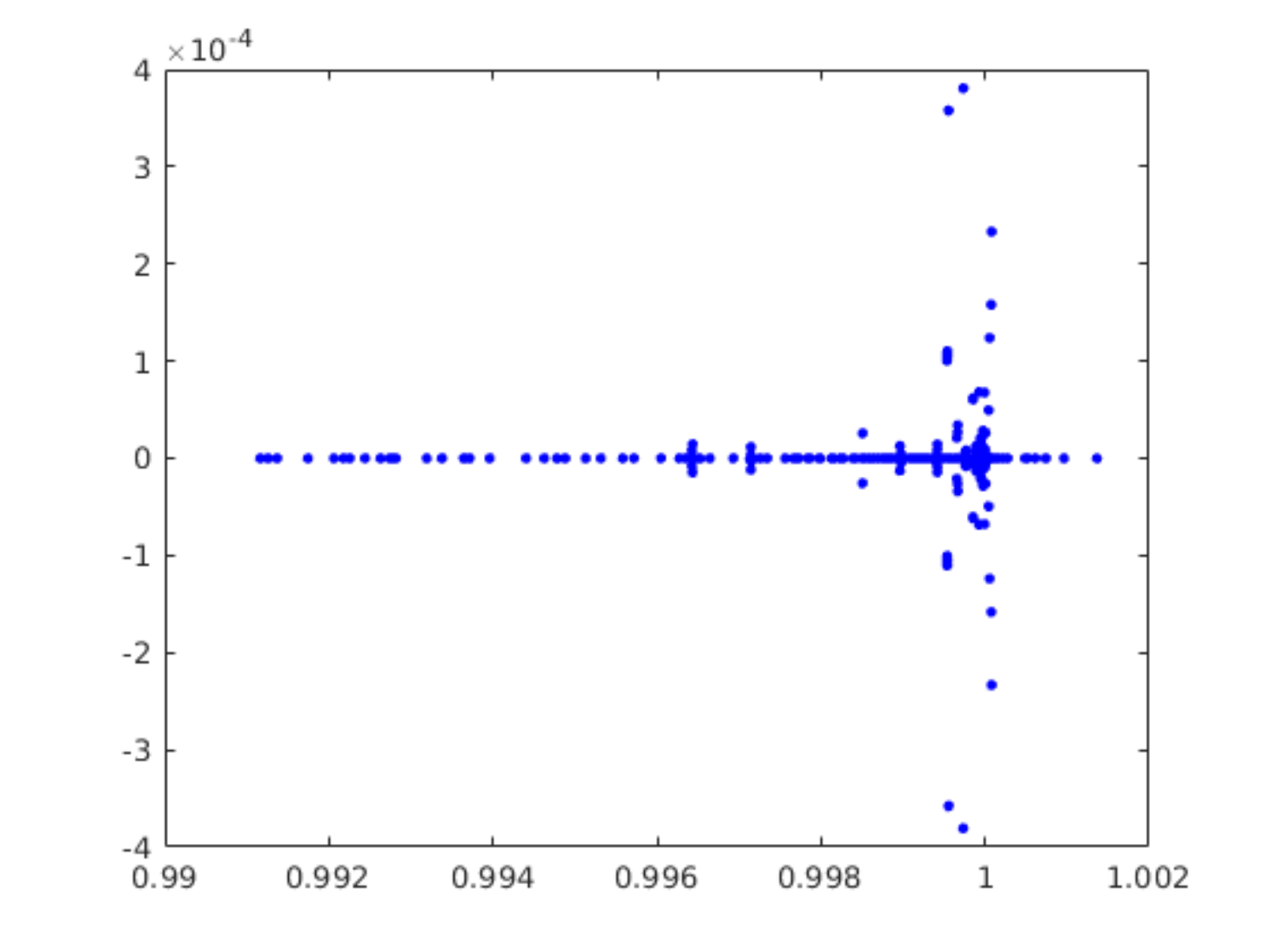}
    \includegraphics[width=0.49\textwidth]{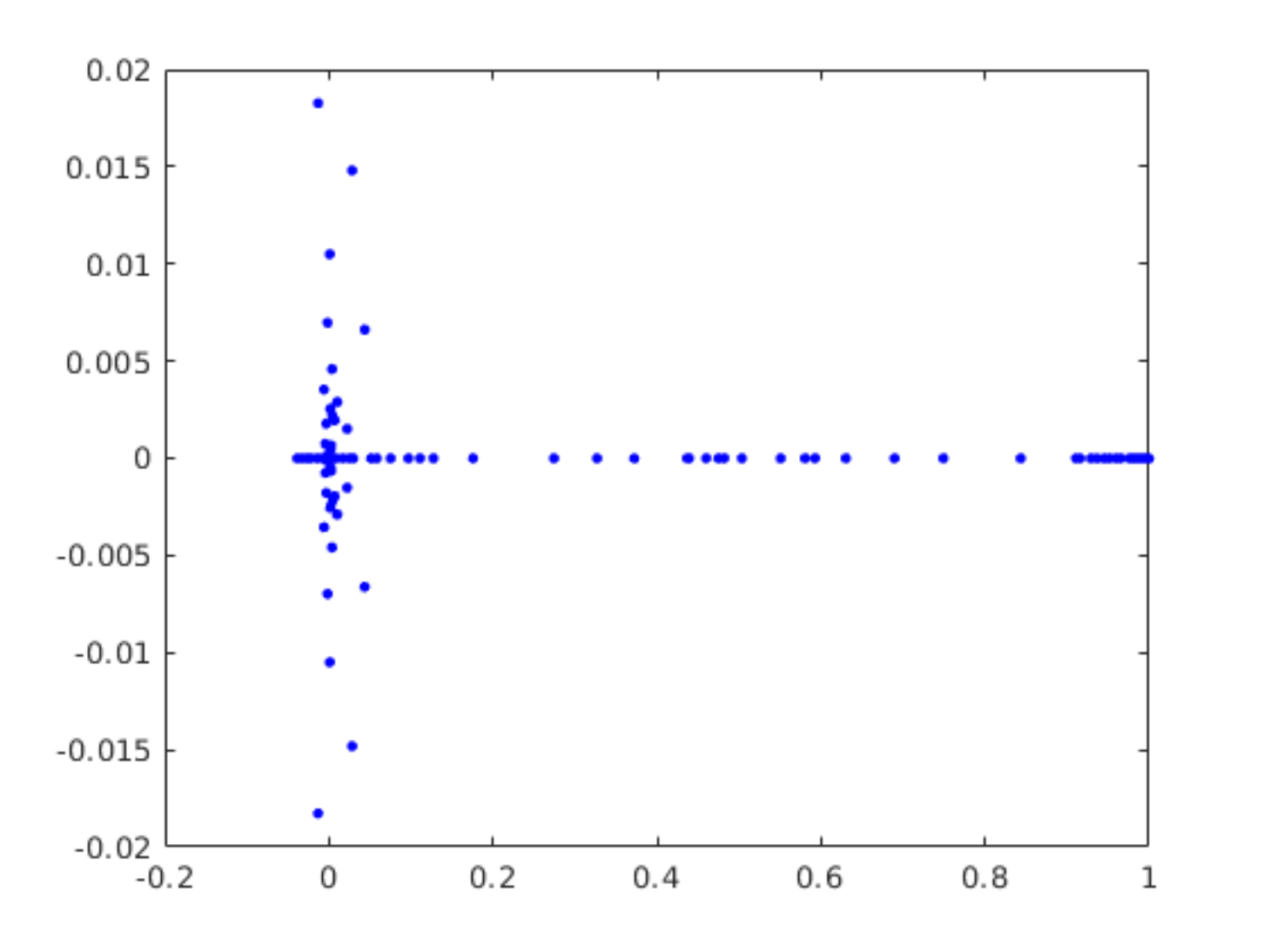}
    \caption{The eigenvalues of the discretization of the heat equation using forward Euler in time with $\Delta t = 10^{-6}$, a grid size of $\Delta x = 0.025$ and $m = 13$ points in the RBF-FD stencil. Left: The eigenvalues using the discretization $I+dtW$, where $I$ is the identity matrix and $W$ is the discretization of the Laplacian. Right: The eigenvalues using the discretization $P+dtW$, where $P$ is the projection matrix and $W$ is the discretization of the Laplacian.}
    \label{figure: eigs}
\end{figure}

Let $\widetilde{U}^n_Z$ be the approximated values of $\tilde u(\cdot,\,t^n)$ on $Z$ and at time $t^n$. This time-stepping scheme is initialized using the initial condition $\widetilde{U}^0_Z:=\tilde u(Z,0) = u(X,0)$.
It is equivalent to consider the discrete equation of the constant-along-normal extended function, and  (\ref{forwardEuler}) becomes
\begin{eqnarray}
   {u}(X,t^{n+1}) &=&  \tilde{u}(X,t^{n})  + \Delta t \Delta \tilde{u}(X,t^n) + \mathcal{O}(\Delta t^2), \quad n\in\mathbb{N}.\label{equation: time discretized PDE}
\end{eqnarray}
We use (\ref{Wu}) and (\ref{Pu}) to approximate the right hand side from the stored approximated solution values $\widetilde{U}^n_Z\approx  \tilde u(Z,t^n)$. Note that we \emph{do not} use pointwise projection, i.e., $  \tilde{u}(X,t^{n}) \approx \widetilde{U}^n_Z$, because $\widetilde{U}^n_Z$ contains discretization and approximation errors. Using (\ref{Pu}) to approximate $\tilde{u}(X,t^{n})$ introduces some averaging into the approximation and increases numerical stability. Indeed, Figure~\ref{figure: eigs} shows the eigenvalues of the discretization of equation~(\ref{equation: time discretized PDE}) on the unit circle for $m=13$ points, a grid size of $\Delta x = 0.025$ and a time step-size of $\Delta t = 10^{-6}$. The scheme that uses pointwise projection leads to an unstable system (eigenvalues larger than $1$) whereas a stable approximation can be achieved using the projection operator.

By design, we have ${u}(X,t^{n+1})=\tilde{u}(Z,t^{n+1})$, whose approximated values will be used to define $\widetilde{U}^{n+1}_Z$.
With all of the above considered, the approximate solution can be updated from time $t^n$ to $t^{n+1}$ by
\begin{equation}\label{RBFCPM}
    \widetilde{U}^{n+1}_Z := (P+\Delta t\, W)  \widetilde{U}^n_Z, \quad n\in\mathbb{N}.
\end{equation}
In (\ref{RBFCPM}), the exterior $Z$ nodes are also used implicitly in the construction of RBF-FD matrices $P$ and $W$ defined in (\ref{Wu}) and (\ref{Pu}) as RBF centers. In other words, the proposed RBF-CPM runs solely based on RBF interpolations with the regularly placed $Z$ nodes.
To evaluate the numerical approximation, say at $X\in\Gamma$ for simplicity, one can evaluate $u(X,t^n)\approx P \widetilde{U}^n_Z$.
Other time-stepping schemes are also possible.
In our experiments, we use the third-order, three-stage SSP Runge-Kutta scheme \cite{shuosher1988}
in advection-dominant problems due to its good linear stability along the imaginary axis.
For simplicity, and to differentiate from other methods, we shall refer to our RBF discretization as the {\it RBF-CPM}.

Given a collection of Cartesian grid points $Z = \{\mathbf{z}_j\}_{j=1}^{n_Z}$ in a small tubular domain $\Omega$ that contains the surface $\Gamma$, the algorithm of the RBF-CPM for a time dependent PDE consists of the following steps:
\begin{enumerate}
    \item Compute the set of surface points $X=\{\mathbf{x}_j\}_{j=1}^{n_Z}\in\Gamma$ via the closest point representation of the surface
 $\Gamma$: $\mathbf{x}_j = cp_\Gamma(\mathbf{z}_j)$, for $\mathbf{z}_j\in Z$, $j=1,\ldots,n_Z$.
	\item Compute the RBF-FD matrices, i.e. the matrices $P$ and $W$ in (\ref{Wu})-(\ref{Pu}). For each surface point $\mathbf{x}_j$, $j=1,2,\ldots,n_Z$:
        \begin{enumerate}
            \item Find the $m$ closest grid points $Z_j = \{\mathbf{z}_{j_i}\}_{i=1}^m$ to $\mathbf{x}_j$.
            \item Compute the RBF-FD weight $\mathbf{w}_j$ at the surface node $\mathbf{x}_j$.
        \end{enumerate}
	\item Solve the surface PDE using an explicit time-stepping scheme, e.g. (\ref{RBFCPM}), until the final time.
\end{enumerate}
Our method allows pre-computation for solving local linear systems; provided that two local neighborhoods $Z_i$ and $Z_j$ share the same geometrical arrangement, the corresponding interpolation matrices are identical, i.e, $A(Z_i,Z_i)=A(Z_j,Z_j)$.
Therefore, an expensive but accurate method can be used.

\subsection{Parameters}
The convergence order of the RBF-FD schemes is limited by the smoothness of the employed kernel and the geometry of the RBF centers \cite{davydov2016optimal}. In this paper, we employ Gaussian RBFs so that the smoothness of the kernel will not be a limiting factor.
To avoid any potential problem of ill-conditioning,
we use the stable RBF-GA method which provides an accurate and stable algorithm and is a cheaper stabilization method over RBF-QR \cite{fornberg2013stable}. These stabilization techniques are independent of the condition number of the matrix $A$ (see Section~\ref{Method description}) using a direct calculation using Gaussian RBFs \cite{fornberg2011stable,fornberg2013stable}.

There are two parameters appearing in the RBF-CPM. These are the shape parameter $\epsilon$ of the Gaussian RBFs and the number of points $m$ in the stencil used locally for the RBF interpolation.
While this section considers the dependence of the numerical method on both parameters, our emphasis will be on the number of points $m$.
The parameter $\epsilon$ was found to have little effect on our results.

We consider two test problems to measure the error of the discrete Laplace-Beltrami operator in comparison to the exact.
The first test problem (P1) approximates the Laplace-Beltrami operator applied to the function $u(\theta) = \sin(\theta)$ on the unit circle.
The relative error in this problem can computed using the known, exact solution $\Delta_\Gamma u = -u$.
In our second problem (P2), the Laplace-Beltrami operator is applied to the function $u(\theta,\phi) = \sin(\phi)$ on the unit sphere $\Gamma$
$$\Gamma = \Big\{\mathbf{x}:\mathbf{x}(\theta,\phi) = (\cos(\theta)\cos(\phi),\sin(\theta)\cos(\phi), \sin(\phi)),-\pi\leq\theta<\pi, -\frac{\pi}{2}\leq\phi\leq\frac{\pi}{2}\Big\}.$$
Here, the relative error can be computed using the known, exact solution $\Delta_\Gamma u = -2u$.

\begin{figure}[h!]
    \centering
    \includegraphics[width=0.49\textwidth]{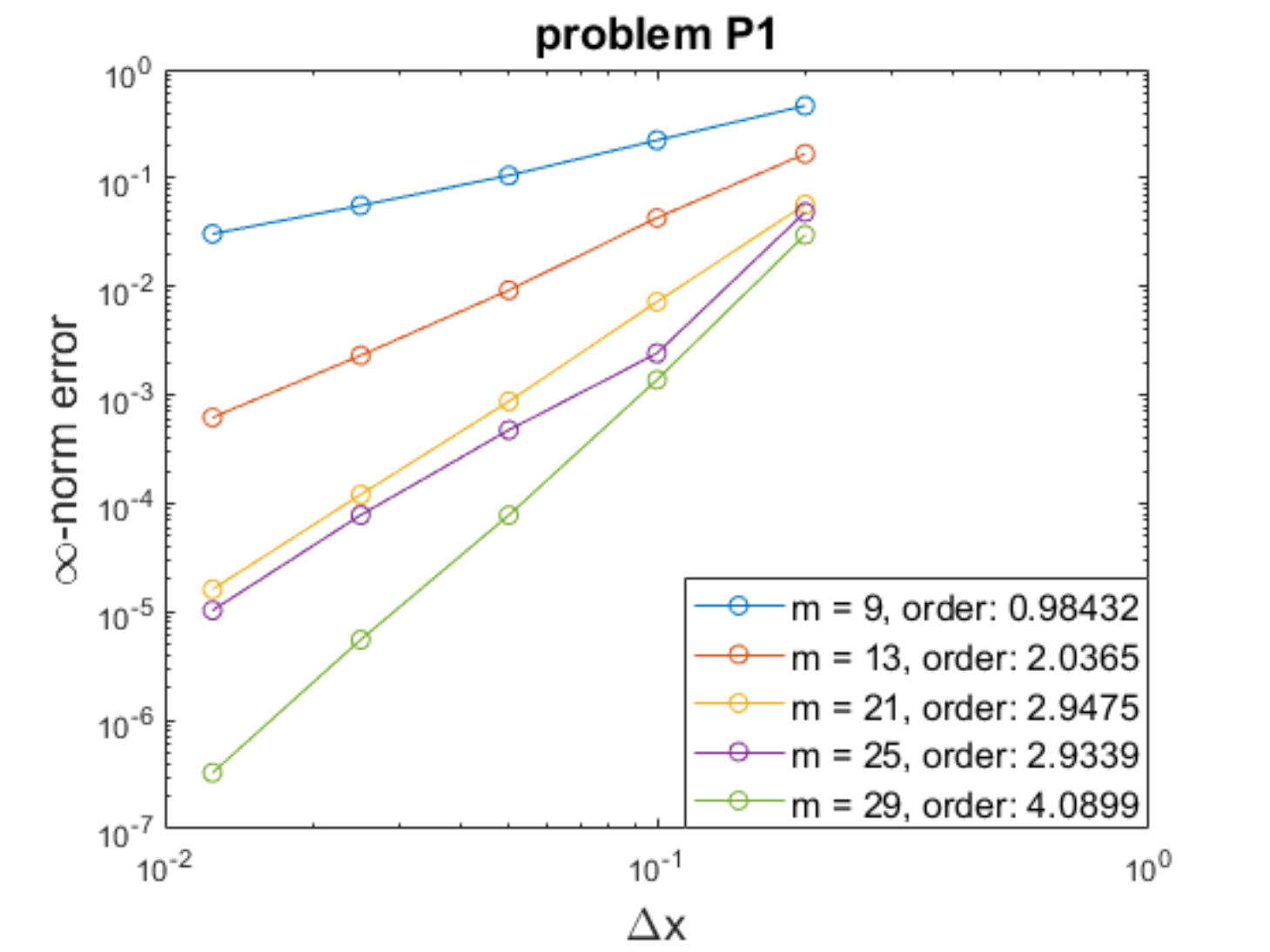}
    \includegraphics[width=0.49\textwidth]{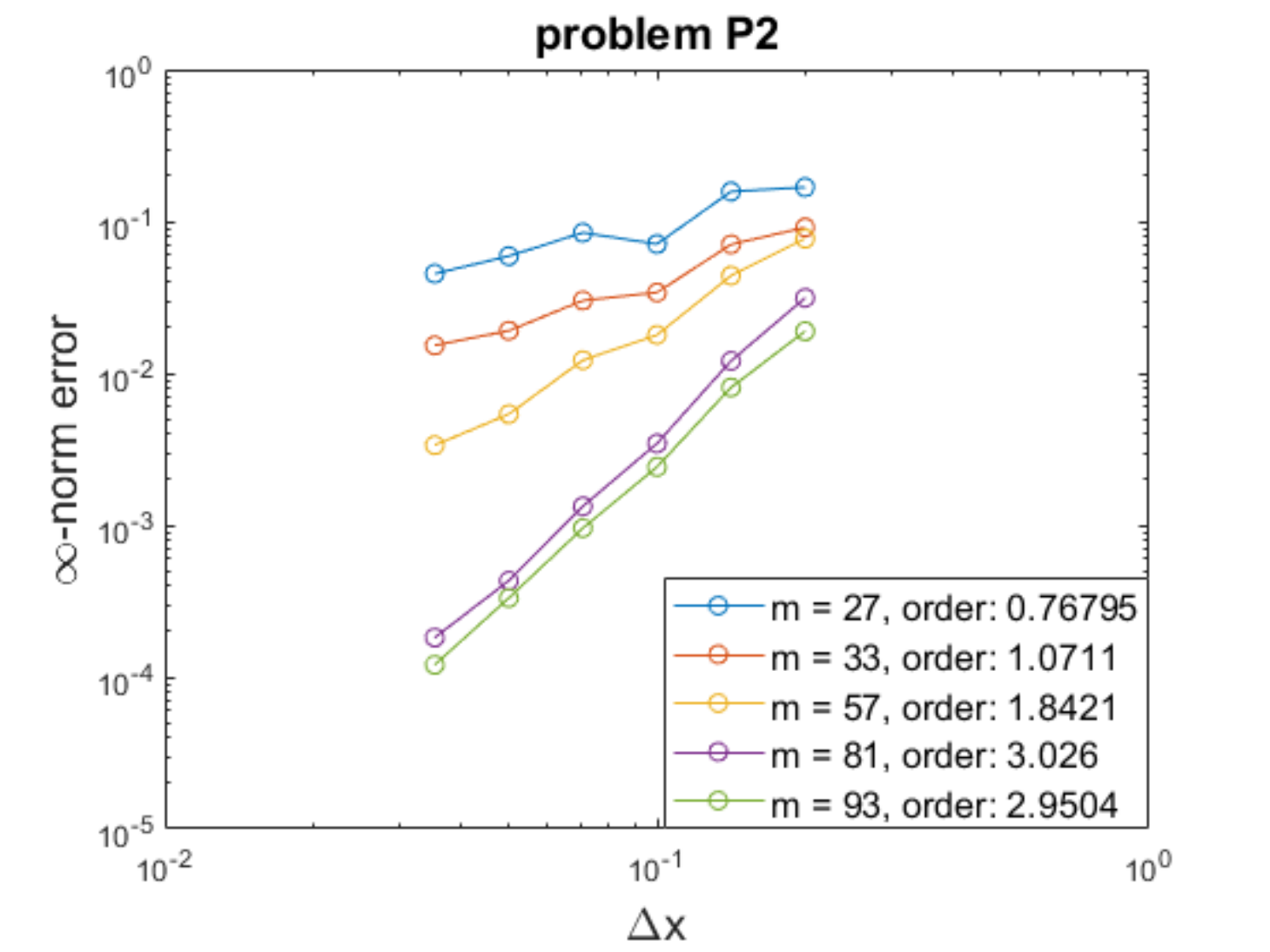}
    \caption{Relative error as a function of grid spacing $\Delta x$ for various stencils
for the approximation of the Laplace-Beltrami operator. Left: problem P1 (two dimensions).  Right: problem P2 (three dimensions). All experiments use $\epsilon=1$.}
    \label{mparameter}
\end{figure}

We apply the RBF-FD method in Equation~(\ref{Wu}) to test problems P1 and P2 for various mesh spacings and selected $m$-values.
See Figure~\ref{mparameter} for the corresponding max norm relative errors.
We find that the observed order of the method increases as the number of points in the finite difference stencil increases.
Due to the regularity of the grid spacing, and the smoothness of the kernel,
our choice of parameter $\epsilon$ has little effect on the results.
In particular, the observed orders of convergence for $\epsilon=1$, $\epsilon=0.1$ and $\epsilon=0.001$ are the same for problems P1 and P2.
For this reason, we simply choose $\epsilon=1$ in the numerical experiments presented in Section~\ref{Numerics}.

\subsection{Computational tube}

The RBF-FD stencils used in the RBF-CPM, i.e., $Z_j\subset Z$ for $j=1,\ldots,n_Z$, are formed
using the $m$-nearest neighboring regular grid points with a predetermined spacing $\Delta x$.
For such problems, the matrix-based formulation of the closest point method \cite{macdonald2009implicit}
can be used to obtain the minimal-sized computational tube.
In this paper, we can take a simpler approach to identify the computational tube by identifying a sufficiently large tube radius $\gamma$.

Consider the Gauss circle problem \cite{weisstein2004gauss}.  In its standard form, it is posed as follows:
Find the number of integer lattice points $m$ inside a circle with radius $r$ centered at the origin.
We shall consider a related formulation:
Find the number of ordered pairs $(x,y)$, with integers $x,y\geq0$, such that
$$x^2+y^2\leq q$$
where the radius of the circle is chosen as $r = \sqrt{q}$, for a fixed integer $q \geq0$. Generalizations to higher dimensions are also available. In three dimensions, the problem uses a sphere centered at the origin.  In this case, we find the number of ordered triplets $(x,y,z)$, with integers $x,y,z\geq0$, such that
$$x^2+y^2+z^2\leq q$$
where the radius of the sphere is $r = \sqrt{q}$, for an integer $q \geq0$.

The integer solutions of the Gauss circle problem and its generalization to three dimensions are given by sequences A057655 (two dimensions) and sequences A117609 (three dimensions) in the On-Line Encyclopedia of Integer Sequences \cite{oeis}. Table~\ref{GaussCircleTable} shows some of the integer solutions for the Gauss circle problem in two and three dimensions.

\begin{table}[h!]
\centering
  \begin{tabular}{|c|c|c|}
    \hline
    $q$ & $m$ (2D) & $m$ (3D)\\ \hline
    0 & 1 & 1\\
    1 & 5 & 7\\
    2 & 9 & 19\\
    3 & 9 & 27\\
    4 & 13 & 33\\
    5 & 21 & 57\\
    \hline
  \end{tabular}
  \caption{The number of lattice points $m$ contained in a ball with radius $r=\sqrt{q}$ in two and three dimensions.}\label{GaussCircleTable}
\end{table}

\begin{figure}[h!]
\centering
    \includegraphics[width=0.3\textwidth]{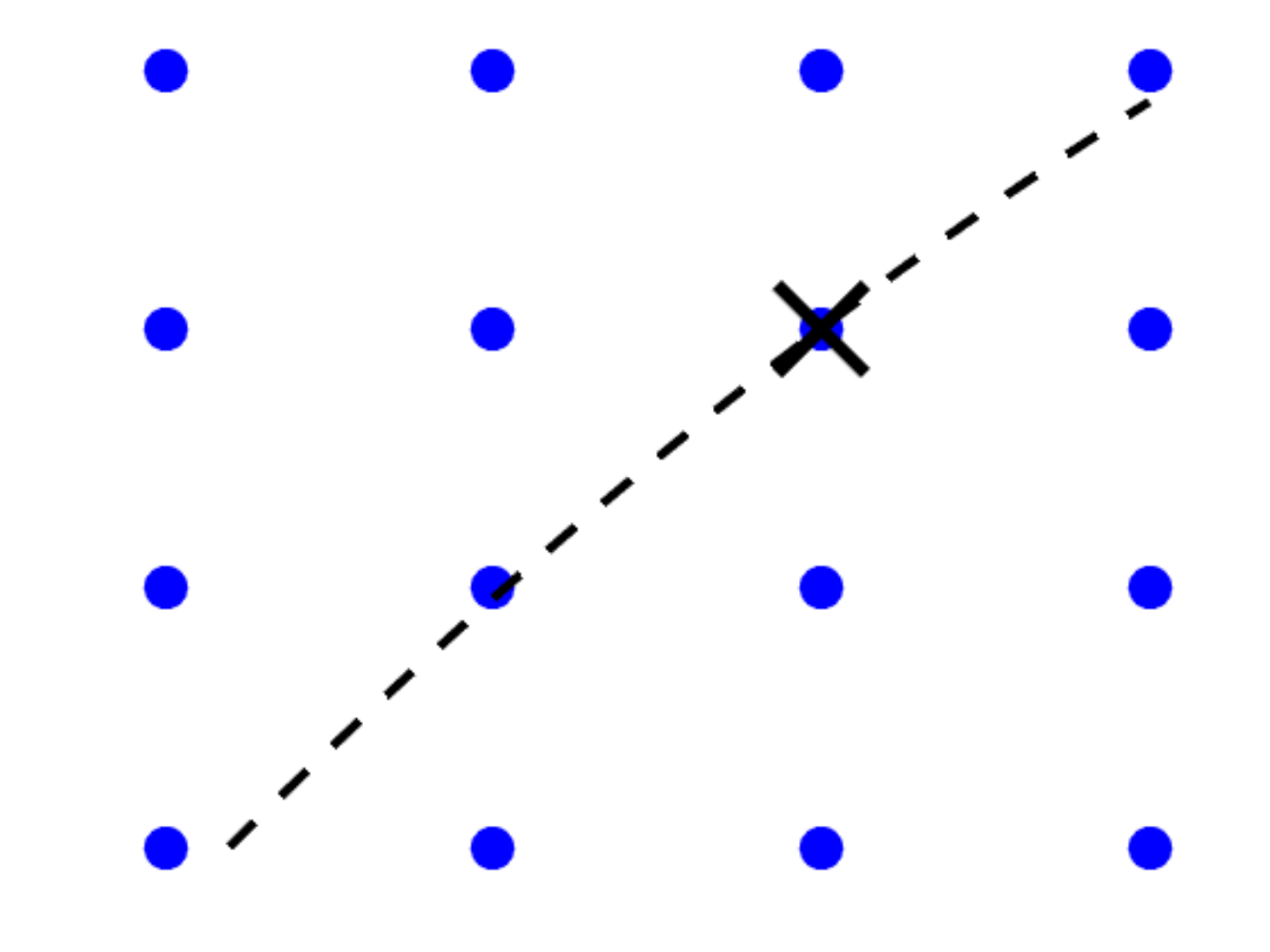}
    \hspace{4cm}
    \includegraphics[width=0.3\textwidth]{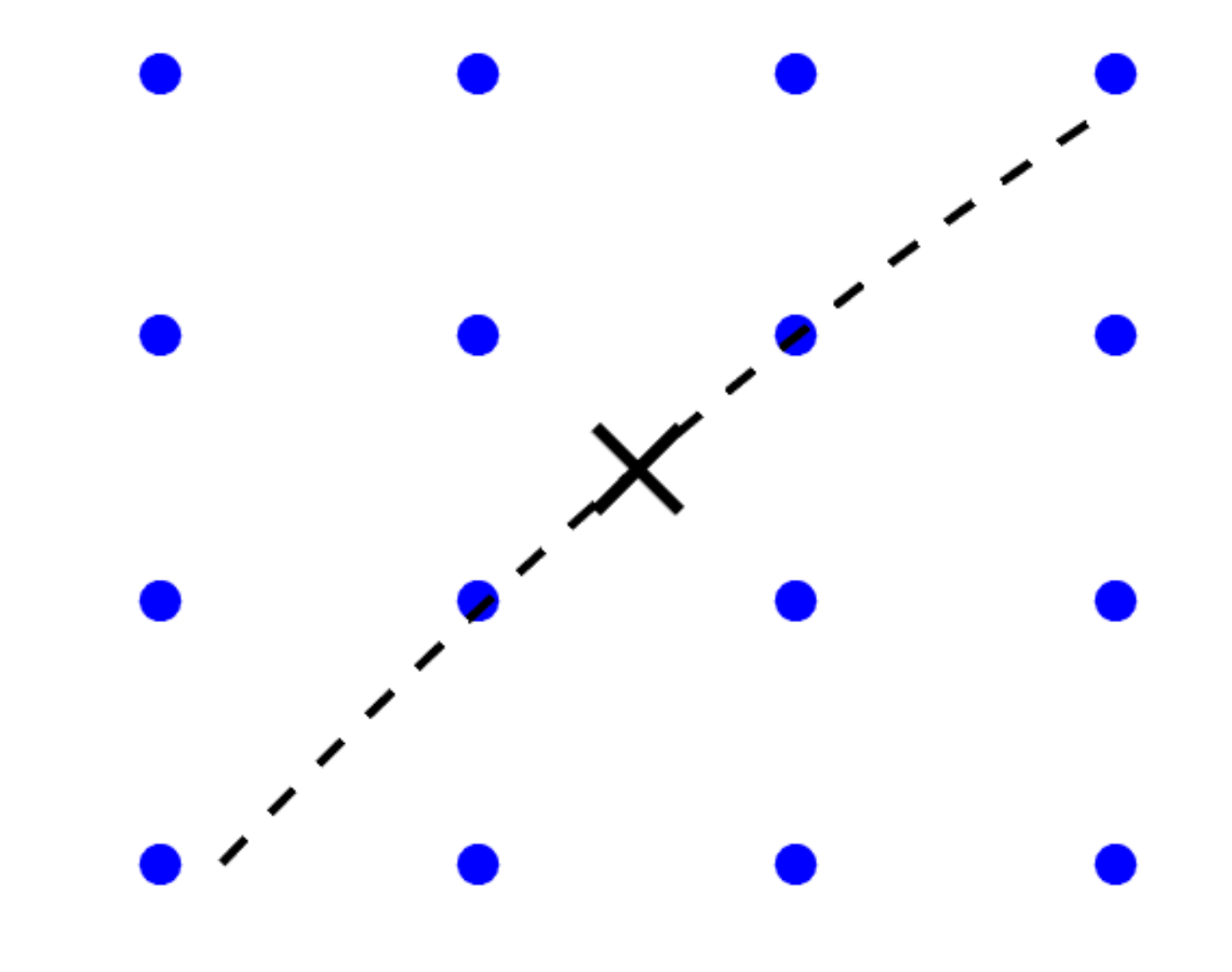}
    \caption{Two cases of a surface point (black x) placement in a mesh grid (blue dots).}
    \label{Cases}
\end{figure}

For a  stencil that uses the $m$ closest grid nodes to a surface point,
a circle can be constructed centered at the surface point that contains these $m$ grid nodes.
In order to construct a computational tube around the surface using the Gauss circle problem,
we need to find a sufficiently large circle independent of the position of the surface point relative to the surrounding grid nodes.
In the optimal case, the surface point lies on a grid node; see Figure~\ref{Cases} (left).
In such an occurrence, the Gauss circle problem can be applied directly, scaled properly with the grid size $\Delta x$.
Otherwise, the surface point does not lie on a grid node.
In order to use the Gauss circle problem, the distance between the surface point and its closest grid node needs to be added to the radius $r$ of the Gauss circle problem.
The worst case appears in Figure~\ref{Cases}, where the surface point lies midway between all four surrounding grid nodes.

To guarantee a sufficiently large circle/sphere radius, we must allow for the worst case.
This leads to a relation between the computational tube radius $\gamma$ (corresponding to the circle/sphere radius $r$) and the number of points $m$ in the RBF-FD stencil;
see Table~\ref{BandwidthTable}.

\begin{table}[h!]
\centering
  \begin{tabular}{|c|c|c|c|}
    \hline
    $\gamma$ (2D) & $m$ (2D) & $\gamma$ (3D) & $m$ (3D)\\ \hline
    ($\sqrt{2}+\sqrt{2}/2)\Delta x$ & 9 & ($\sqrt{3}+\sqrt{3}/2)\Delta x$ & 27\\
    (2+$\sqrt{2}/2)\Delta x$ & 13 & (2+$\sqrt{3}/2)\Delta x$ & 33\\
    $(\sqrt{5}+\sqrt{2}/2)\Delta x$ & 21 & $(\sqrt{5}+\sqrt{3}/2)\Delta x$ & 57\\
    $(\sqrt{8}+\sqrt{2}/2)\Delta x$ & 25 & $(\sqrt{6}+\sqrt{3}/2)\Delta x$ & 81\\
    & & $(\sqrt{8}+\sqrt{3}/2)\Delta x$ & 93\\
    \hline
  \end{tabular}
  \caption{Computational tube radius $\gamma$ for an $m$-point RBF-FD stencil in two and three dimensions.}\label{BandwidthTable}
\end{table}

\section{Numerical experiments}\label{Numerics}
In this section, we test our method, the RBF-CPM,  on a number of examples in two and three dimensions. For problems involving the Laplace-Beltrami operator, we choose $m=13$ in two dimensions and $m=57$ in three dimensions, which provides a second order approximation for smooth solutions (cf. Figure~\ref{mparameter}). The radius of the computational tube is set according to the values specified in Table~\ref{BandwidthTable}. Finding the optimal RBF shape parameter is out of the scope of this paper, and thus we set $\epsilon=1$ (unscaled Gaussian RBFs are used). Unless stated otherwise, forward Euler with a time step-size $\Delta t=0.1\Delta x^2$ is used for the discretization of the time derivatives.

\subsection{Examples in two dimensions}
First, we perform numerical experiments for applications in two dimensions to test the convergence of the proposed method.

\subsubsection{Heat equation on a circle}
In our first experiment, we consider the heat equation
$$u_t=\Delta_\Gamma u$$
intrinsic to the unit circle $\Gamma$. Following \cite{ruuth2008simple}, for an initial profile $u(\theta,0)=\sin\theta$, the exact solution is
$$u(\theta,t)=e^{-t}\sin\theta.$$
Using an analytic closest point representation of the unit circle, the surface heat equation is discretized and solved using Equation~(\ref{RBFCPM}).
Table~\ref{errorCircleTableExplicit} shows the relative errors as well as the convergence rates for different grid sizes $\Delta x$ and number of points $N$ on the computational domain.
The convergence rate here appears to be at least second-order. Using the original closest point method, the number of points
in the computational tube that are required for a second order approximation with $\Delta x = 0.00625$ is $7276$, whereas our method uses $5464$ points.  This corresponds to a reduction of $25\%$.

\begin{table}[h!]
\centering
  \begin{tabular}{|l|c|c|c|}
    \hline
    $\Delta x$ & $N$ & Rel. error ($t = 1$) & Conv. rates \\ \hline
    0.2 & 172 & 7.15$\times 10^{-3}$ & - \\
    0.1 & 336 & 1.22$\times 10^{-3}$ & 2.55  \\
    0.05 & 688 & 2.23$\times 10^{-4}$ & 2.46  \\
    0.025 & 1376 & 5.15$\times 10^{-5}$ & 2.11  \\
    0.0125 & 2708 & 1.35$\times 10^{-5}$ & 1.93  \\
    0.00625 & 5464 & 3.15$\times 10^{-6}$ & 2.10  \\
    \hline
  \end{tabular}
  \caption{Relative errors and convergences rates for the approximate solution at time $t=1$ for the heat equation on the unit circle. Errors are measured in the infinity norm.}
  \label{errorCircleTableExplicit}
\end{table}

\subsubsection{Heat equation on a semicircle}
\begin{figure}
    \centering
    \includegraphics[width=0.5\textwidth]{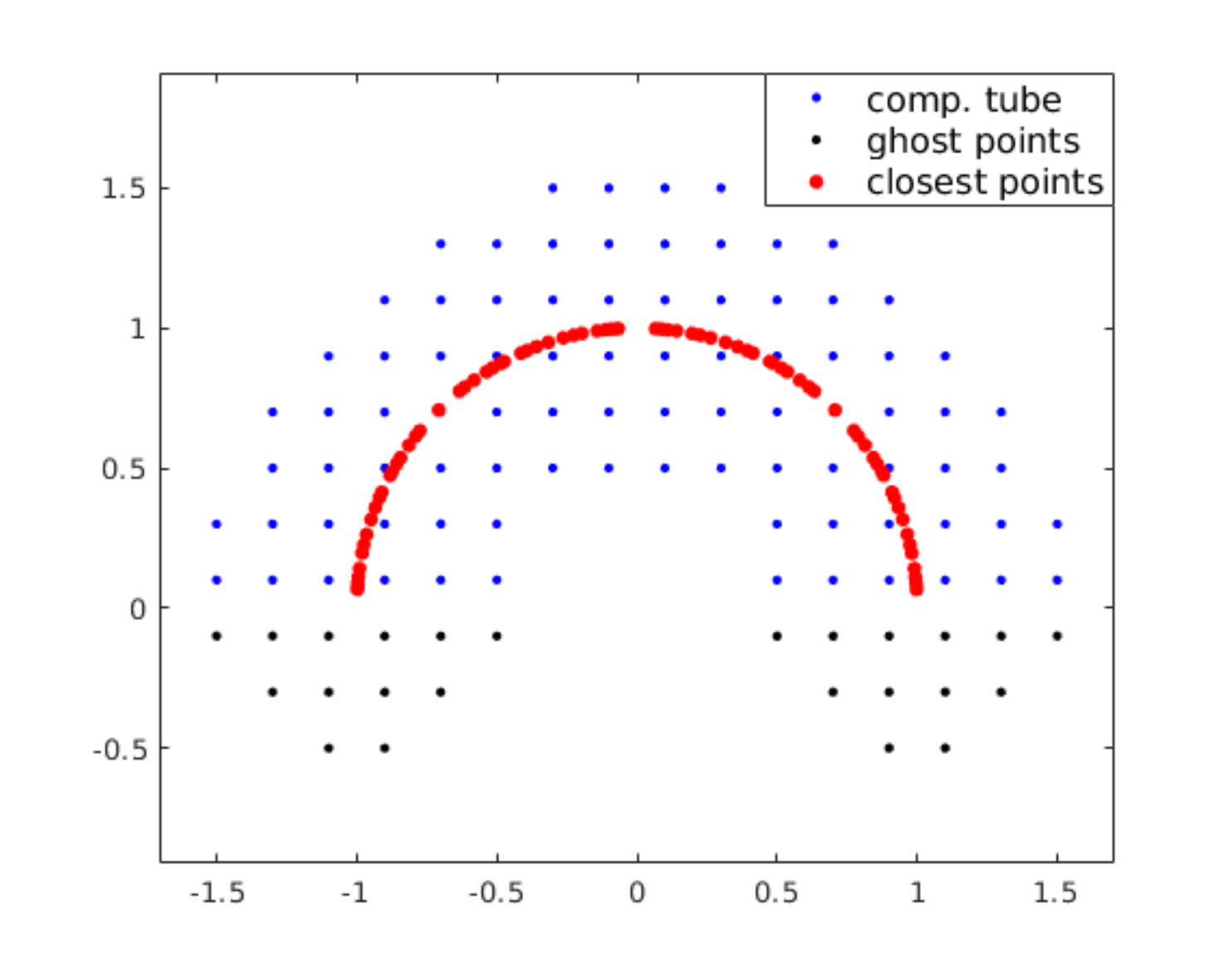}
    \caption{The computational tube around the semicircle, including internal grid points (blue dots) and ghost points (black dots), along with their corresponding closest points on the semicircle (red dots).}
    \label{figure: heatSemicircleCompTube}
\end{figure}

Next, we consider the surface heat equation
on the unit semicircle with homogeneous Dirichlet boundary conditions. Given an initial profile $u(\theta,0)=\sin\theta$, the exact solution is
$$u(\theta,t)=e^{-t}\sin\theta.$$
Following \cite{macdonald2011solving}, we introduce a modified closest point mapping
$\overline{cp}_\Gamma(\mathbf{z}) = cp_\Gamma(2cp_\Gamma(\mathbf{z})-\mathbf{z})$ which equals $cp_\Gamma(\mathbf{z})$ for points $\mathbf{z}$ that map to the interior of the semi-circle.
Grid nodes that satisfy $\overline{cp}_\Gamma(\mathbf{z})\neq cp_\Gamma(\mathbf{z})$ are
called ghost points $\mathbf{z}_g$.  At such points, the function $u$ is extended by $-u(\overline{cp}_\Gamma(\mathbf{z}_g))$. Figure~\ref{figure: heatSemicircleCompTube} shows the computational tube used in this example.

Table~\ref{errorSemicircleTableExplicit} presents
the relative errors as well as the convergence rates for different grid sizes $\Delta x$ and number of points $N$ on the computational domain.
Second-order convergence is observed.
Using the original closest point method, the number of points in the computational tube that are required for a second order approximation at $\Delta x = 0.00625$ is $3682$,
whereas the RBF-CPM uses $2756$ points.

\begin{table}[h!]
\centering
  \begin{tabular}{|l|c|c|c|}
    \hline
    $\Delta x$ & $N$ & Rel. error ($t = 1$) & Conv. rates \\ \hline
    0.2 & 110 & 7.38$\times 10^{-3}$ & - \\
    0.1 & 192 & 1.14$\times 10^{-3}$ & 2.69  \\
    0.05 & 368 & 2.12$\times 10^{-4}$ & 2.43  \\
    0.025 & 712 & 5.02$\times 10^{-5}$ & 2.08  \\
    0.0125 & 1378 & 1.34$\times 10^{-5}$ & 1.91  \\
    0.00625 & 2756 & 3.13$\times 10^{-6}$ & 2.09  \\
    \hline
  \end{tabular}
  \caption{Relative errors and convergence rates for the approximate solution at time $t=1$ for the heat equation on the unit semicircle.  Errors are measured in the infinity norm.}
  \label{errorSemicircleTableExplicit}
\end{table}

\subsubsection{Advection equation on an ellipse}
The next example is the advection equation on an ellipse. Following \cite{ruuth2008simple}, the equation
$$u_t+u_s=0,$$
with $s$ being the arclength, is imposed on an ellipse with major axis $b=1.25$ along the y-axis and minor axis $a=0.75$ along the x-axis. By \cite{ruuth2008simple}, application of the closest point principles to the surface PDE leads to the embedding PDE
$$u_t+\mathbf{T}(x,y)\cdot\nabla u=0$$
with
$$\mathbf{T}(x,y)=\frac{(-y/b^2,x/a^2)}{\sqrt{y^2/b^4+x^2/a^4}}.$$
For an initial profile $u(s,0)=\sin^3(2\pi s/L)$, the exact solution for subsequent times $t\ge 0$ is
$$u(s,t) = \sin^3(2\pi (s-t)/L),$$
where $L$ is the length of the perimeter of the ellipse. Using a parametrization for the ellipse, the closest point representation is calculated via optimization techniques. Due to the generic centered nature of the RBF-FD stencils used for approximating the first-order derivatives, the TVD-RK3 scheme \cite{gottlieb1998total} is chosen for the time discretization with a time step-size $\Delta t=0.5\Delta x$. In this example, the computational tube radius is $\gamma=(\sqrt{2}+\sqrt{2}/2)\Delta x$ with $m=9$ points in the finite difference stencil. Table~\ref{errorAdvEllipseTable} shows the error at the final time $t=1$ and the estimated order of convergence of the method for various grid spacings and number of points in the computational domain.
Second-order convergence is observed. The number of points in the computational tube that are required for a second order approximation at $\Delta x = 0.00625$ using the original closest point method is $7276$.   Our method uses $41\%$ fewer points.

\begin{table}[h!]
\centering
  \begin{tabular}{|l|c|c|c|}
    \hline
    $\Delta x$ & $N$ & Rel. error ($t = 1$) & Conv. rates \\ \hline
    0.2 & 136 & 8.99$\times 10^{-2}$ & - \\
    0.1 & 272 & 9.80$\times 10^{-3}$ & 3.20  \\
    0.05 & 552 & 2.25$\times 10^{-3}$ & 2.12  \\
    0.025 & 1080 & 5.59$\times 10^{-4}$ & 2.01  \\
    0.0125 & 2168 & 1.40$\times 10^{-4}$ & 2.00  \\
    0.00625 & 4332 & 3.51$\times 10^{-5}$ & 1.99  \\
    0.003125 & 8652 & 8.75$\times 10^{-6}$ & 2.00 \\
    \hline
  \end{tabular}
  \caption{Relative errors and convergences rates for the approximate solution at time $t=1$ for the advection equation on an ellipse.  Errors are measured in the infinity norm.}
  \label{errorAdvEllipseTable}
\end{table}

\subsubsection{Advection-diffusion equation on an ellipse}
Our next example considers advection-diffusion on an ellipse. The equation
$$u_t+u_s=u_{ss},$$
with $s$ being the arclength, is imposed on an ellipse with major axis $b=1.25$ (along the y-axis) and minor axis $a=0.75$ (along the x-axis).
Similar to the previous example, application of the closest point principles to the surface PDE gives
$$u_t+\mathbf{T}(x,y)\cdot\nabla u=\Delta u$$
with
$$\mathbf{T}(x,y)=\frac{(-y/b^2,x/a^2)}{\sqrt{y^2/b^4+x^2/a^4}}.$$
For an initial profile $u(s,0)=\sin(2\pi s/L)$, the exact solution has the form
$$u(s,t)=e^{-2\pi t/L}\sin(2\pi(s-t)/L),$$
where $L$ is the length of the perimeter of the ellipse.
Table~\ref{errorAdvDifEllipseTable} shows the relative error of the approximate solution compared to the exact as well as the estimated order of convergence. Second-order convergence is observed.

\begin{table}[h!]
\centering
  \begin{tabular}{|l|c|c|c|}
    \hline
    $\Delta x$ & $N$ & Rel. error ($t = 1$) & Conv. rates \\ \hline
    0.2 & 172 & 9.66$\times 10^{-3}$ & - \\
    0.1 & 348 & 1.34$\times 10^{-3}$ & 2.85  \\
    0.05 & 692 & 4.88$\times 10^{-4}$ & 1.46  \\
    0.025 & 1384 & 1.25$\times 10^{-4}$ & 1.97  \\
    0.0125 & 2792 & 2.72$\times 10^{-5}$ & 2.20  \\
    0.00625 & 5552 & 6.86$\times 10^{-6}$ & 1.99  \\
    0.003125 & 11100 & 1.68$\times 10^{-6}$ & 2.03  \\
    \hline
  \end{tabular}
  \caption{Relative errors and convergence rates for the approximate solution at time $t=1$ for the advection-diffusion equation on an ellipse.  Errors are measured in the infinity norm.}
  \label{errorAdvDifEllipseTable}
\end{table}

\subsection{Examples in three dimensions}
Next, we apply our proposed method to examples in three dimensions.

\subsubsection{Heat equation on a sphere}
For our first three dimensional example, consider the heat equation
$$u_t=\Delta_\Gamma u$$
on the unit sphere $\Gamma$. For an initial profile $u(\theta,\phi,0)=\sin\phi$, the exact solution for all times $t$ is
$$u(\theta,\phi,t) = e^{-2t}\sin\phi.$$
An analytic closest point representation of the unit sphere is used and the surface heat equation is discretized and solved using Equation~(\ref{RBFCPM}).
Table~\ref{errorSphereTableExplicit} shows the relative errors as well as the convergence rates for different grid sizes $\Delta x$ and number of points $N$ in the computational tube. Second-order convergence is observed. Using $\Delta x = 0.0125$ and a second order finite difference discretization in the original closest point method leads to
$663880$ points in the computational domain, while the RBF-CPM uses $498392$ points.  This corresponds to a reduction of $25\%$.

\begin{table}[h!]
\centering
  \begin{tabular}{|l|c|c|c|}
    \hline
    $\Delta x$ & $N$ & Rel. error ($t = 1$) & Conv. rates \\ \hline
    0.2 & 2240 & 8.18$\times 10^{-3}$ & - \\
    0.1 & 8072 & 2.21$\times 10^{-3}$ & 1.89  \\
    0.05 & 31416 & 5.42$\times 10^{-4}$ & 2.03  \\
    0.025 & 125216 & 1.36$\times 10^{-4}$ & 1.99  \\
    0.0125 & 498392 & 3.40$\times 10^{-5}$ & 2.00  \\
    \hline
  \end{tabular}
  \caption{Relative errors and convergence rates for the approximate solution at time $t=1$ for the heat equation on the unit sphere.  Errors are measured in the infinity norm.}
  \label{errorSphereTableExplicit}
\end{table}

\subsubsection{Advection on a torus}
In this example, the solution of the advection equation on a torus is approximated. Following \cite{greer2006improvement}, for a torus defined as
\small
$$\Gamma = \Big\{\mathbf{x}:\mathbf{x}(\phi,\theta) = \left(\Big(\frac{1}{2}\cos(\phi) + 1\Big)\cos(\theta),\Big(\frac{1}{2}\cos(\phi) + 1\Big)\sin(\theta),\frac{1}{2}\sin(\phi)\right),-\pi\leq\theta,\phi\leq\pi\Big\},$$
\normalfont
the advection equation is given by
$$u_t+u_\phi=0.$$
For an initial profile
$$u(\phi,\theta,0) = f(\phi) = \left\{
                \begin{array}{ll}
                  g(\frac{\phi+\pi}{\pi}), & \hbox{$-\pi\leq\phi\leq0$,} \\
                  g(\frac{\pi-\phi}{\pi}), & \hbox{$0<\phi<\pi$,}
                \end{array}
              \right.$$
where
$$g(s) = \frac{e^{1/(s-1)}-e^{-1/s}}{e^{-1/s}+e^{1/(s-1)}}, $$
the exact solution at time $t$ is
$$u(\phi,\theta,t) = f(\phi-t).$$

\begin{table}[h!]
\centering
  \begin{tabular}{|l|c|c|c|}
    \hline
    $\Delta x$ & $N$ & Rel. error ($t = 1$) & Conv. rates \\ \hline
    0.1 & 11392 & 1.76$\times 10^{-2}$ & -  \\
    0.05 & 45464 & 2.99$\times 10^{-3}$ & 2.56  \\
    0.025 & 181480 & 4.88$\times 10^{-4}$ & 2.62  \\
    0.0125 & 725200 & 9.52$\times 10^{-5}$ & 2.36  \\
    0.00625 & 2901248 & 1.80$\times 10^{-5}$ & 2.40  \\
    \hline
  \end{tabular}
  \caption{Relative error and  convergence rates for the approximate solution at time $t=1$ for the advection equation on a torus.  Errors are measured in the infinity norm.}
  \label{errorTorusTable}
\end{table}

In this example, an analytic closest point representation of the torus is used.
The computational tube radius is $\gamma=(2+\sqrt{3}/2)\Delta x$ and a stencil with $m=33$ points is chosen.
The stable TVD-RK3 scheme is chosen for the discretization in time with step-size $\Delta t = 0.5\Delta x$. Table~\ref{errorTorusTable} shows the error and the convergence rate of the method for various grid sizes $\Delta x$.
Convergence is at least second-order. The number of points in the computational tube that are required for a second order approximation at $\Delta x = 0.00625$ using the original closest point method is $4163904$.  Our method uses $30\%$ fewer points.

\subsubsection{Image denoising on a sphere}
Our next example concerns image denoising for a textured image on the unit sphere.
Following \cite{biddle2013volume}, we apply the Perona-Malik model to denoise surface images.
The equation has the form
$$u_t = \nabla_\Gamma\cdot (g(|\nabla_\Gamma u|)\nabla_\Gamma u)$$
where $g$ is the diffusion coefficient given by
$$g(s) = \frac{1}{1+(s/\lambda)^2}$$
and $\lambda$ is a coefficient. The parameter $\lambda$ and the final time $t$ of the computation control the denoising of an image.

\begin{figure}
    \centering
    \includegraphics[width=0.5\textwidth]{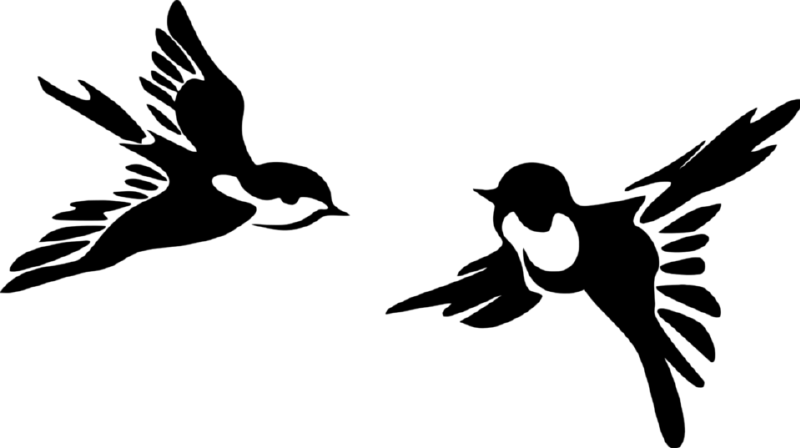}
    \includegraphics[width=\textwidth]{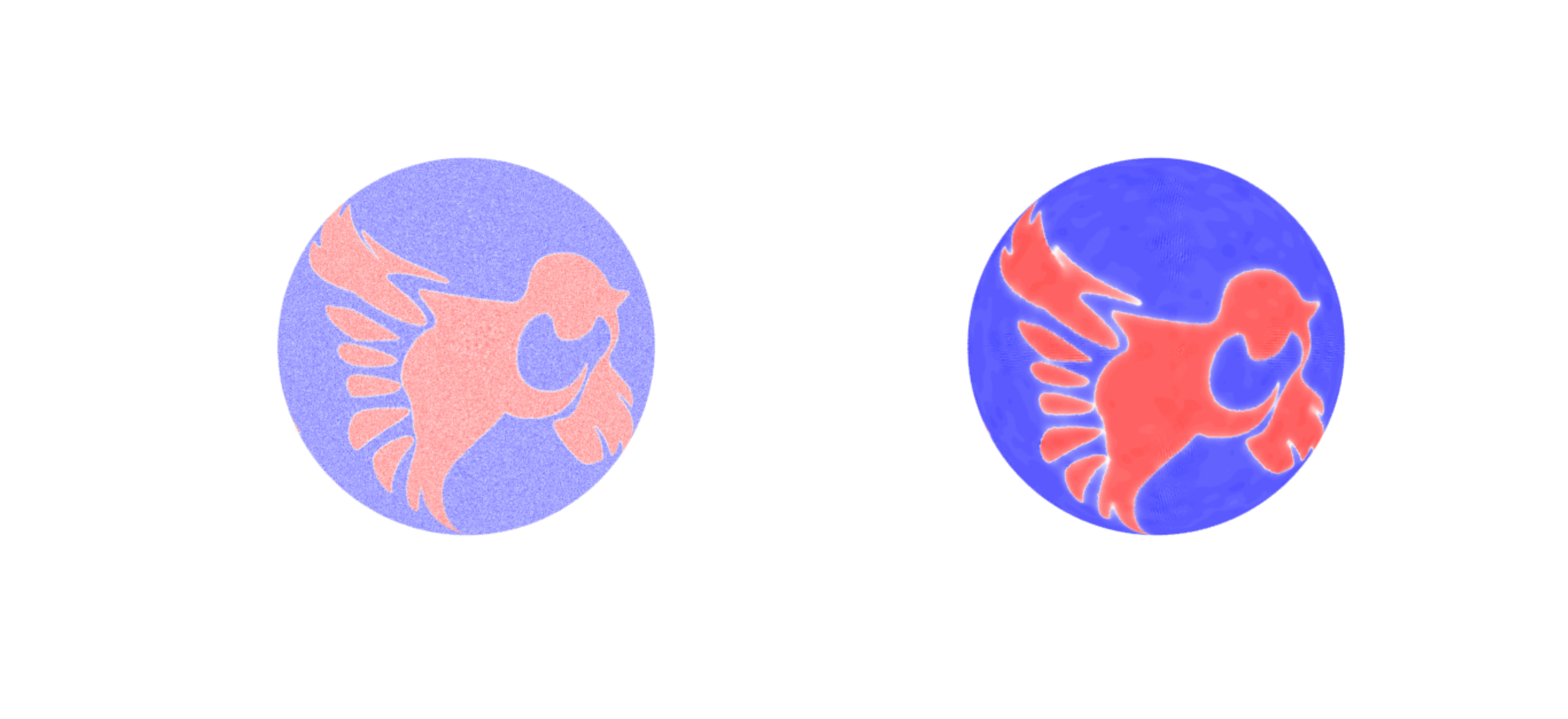}
    \includegraphics[width=\textwidth]{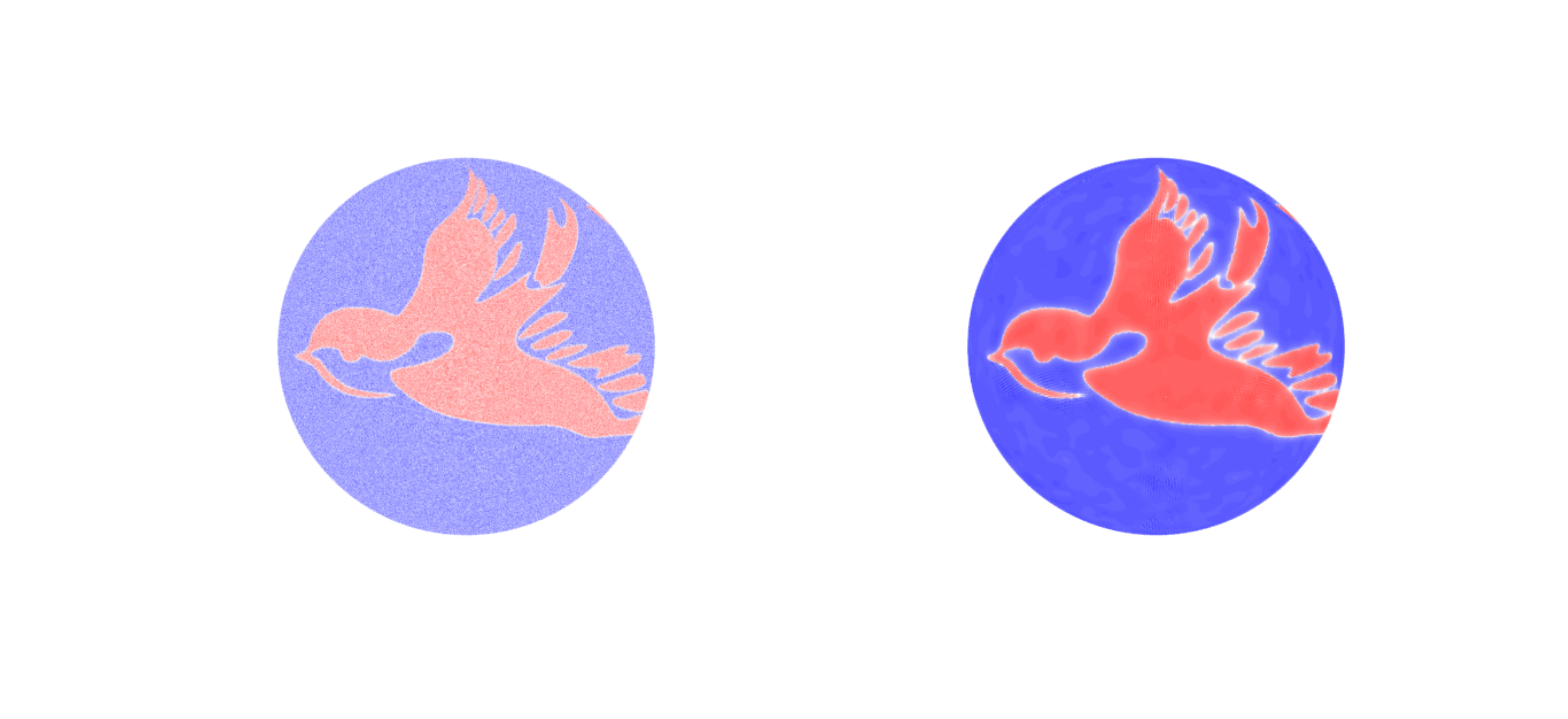}
    \caption{The initial image (top) is warped onto a sphere and Gaussian noise is added. The noisy image (left) and the denoised image (right) are shown after 120 iterations on 1280598 points.}
    \label{denoising}
\end{figure}

To construct the initial image, we add Gaussian noise (zero mean and 0.2 standard deviation) to the image of two birds.
The noisy image is scaled to the interval $[0,1]$. Figure~\ref{denoising} shows the results after denoising the image using $\lambda = 5$ and  $120$ time steps with $\Delta t = 0.2\Delta x^2$
and $\Delta x=0.005$. We find this choice for the parameter $\lambda$ is sufficiently small to preserve edges in the denoised image.
In this example, the computational tube radius is chosen to be $\gamma=(2+\sqrt{3}/2)\Delta x$ with $m=33$ points in the stencil.  The number of points in the image is $1280598$.

\subsubsection{Reaction-diffusion systems}
Our final example evolves the Gray-Scott reaction-diffusion model \cite{GRAY19841087} on a triangulated surface.    The Gray Scott model describes the chemical reaction
$$U+2V\longrightarrow3V, $$
$$V\longrightarrow P,$$
where $U$, $V$ and $P$ are chemicals.
The corresponding surface model has the form
$$u_t = F(1-u)-uv^2+\mathcal{D}_u\Delta_Su$$
$$v_t = -(F+k)v+uv^2+\mathcal{D}_v\Delta_Sv$$
where $u$, $v$ are the concentrations of the chemicals, $\mathcal{D}_u$, $\mathcal{D}_v$ are the diffusion rates, $k$ is the conversion rate from $V$ to $P$ and $F$ is the feed rate of $U$.

\begin{figure}
    \centering
    \includegraphics[width=0.49\textwidth]{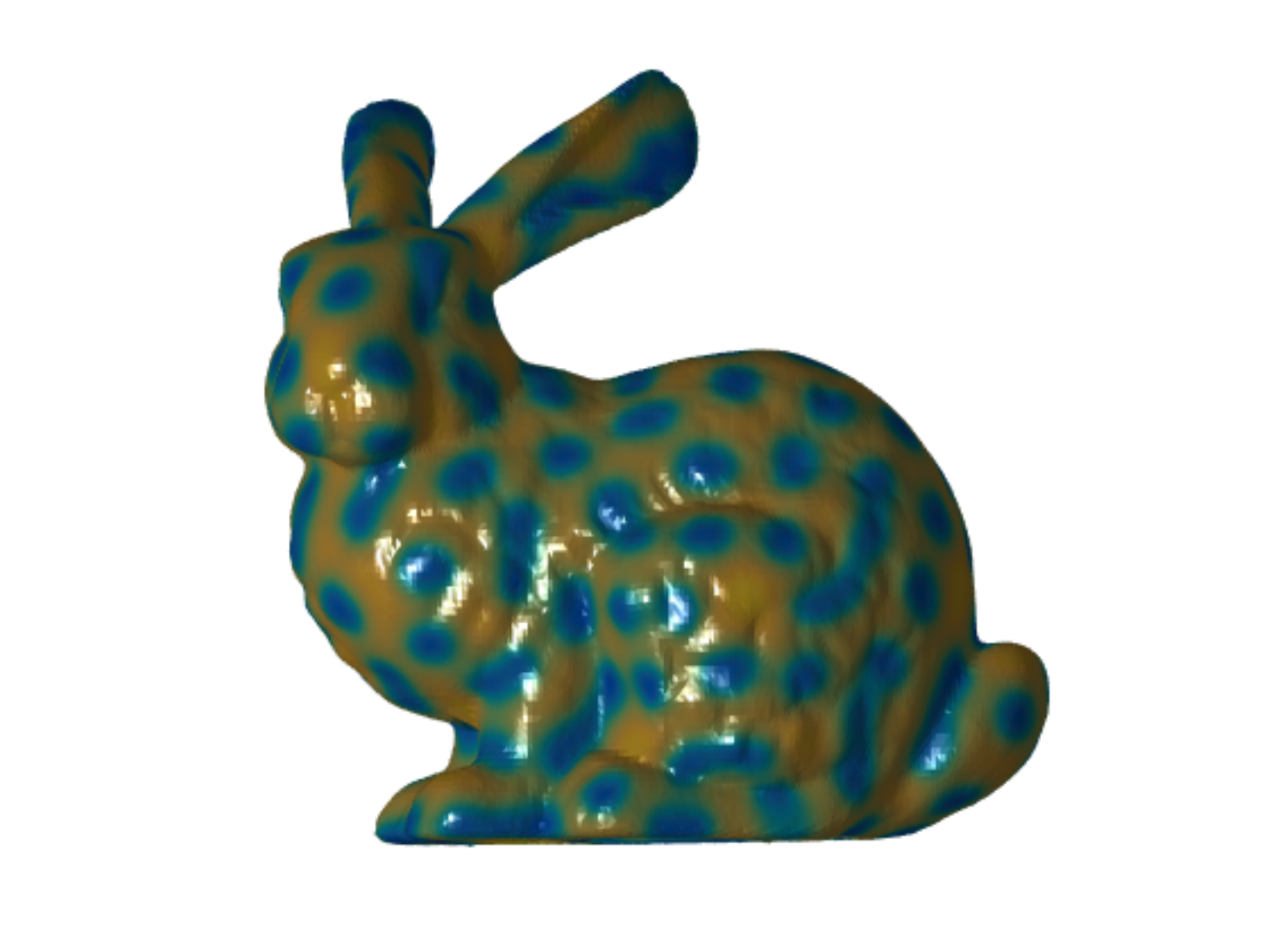}
    \includegraphics[width=0.49\textwidth]{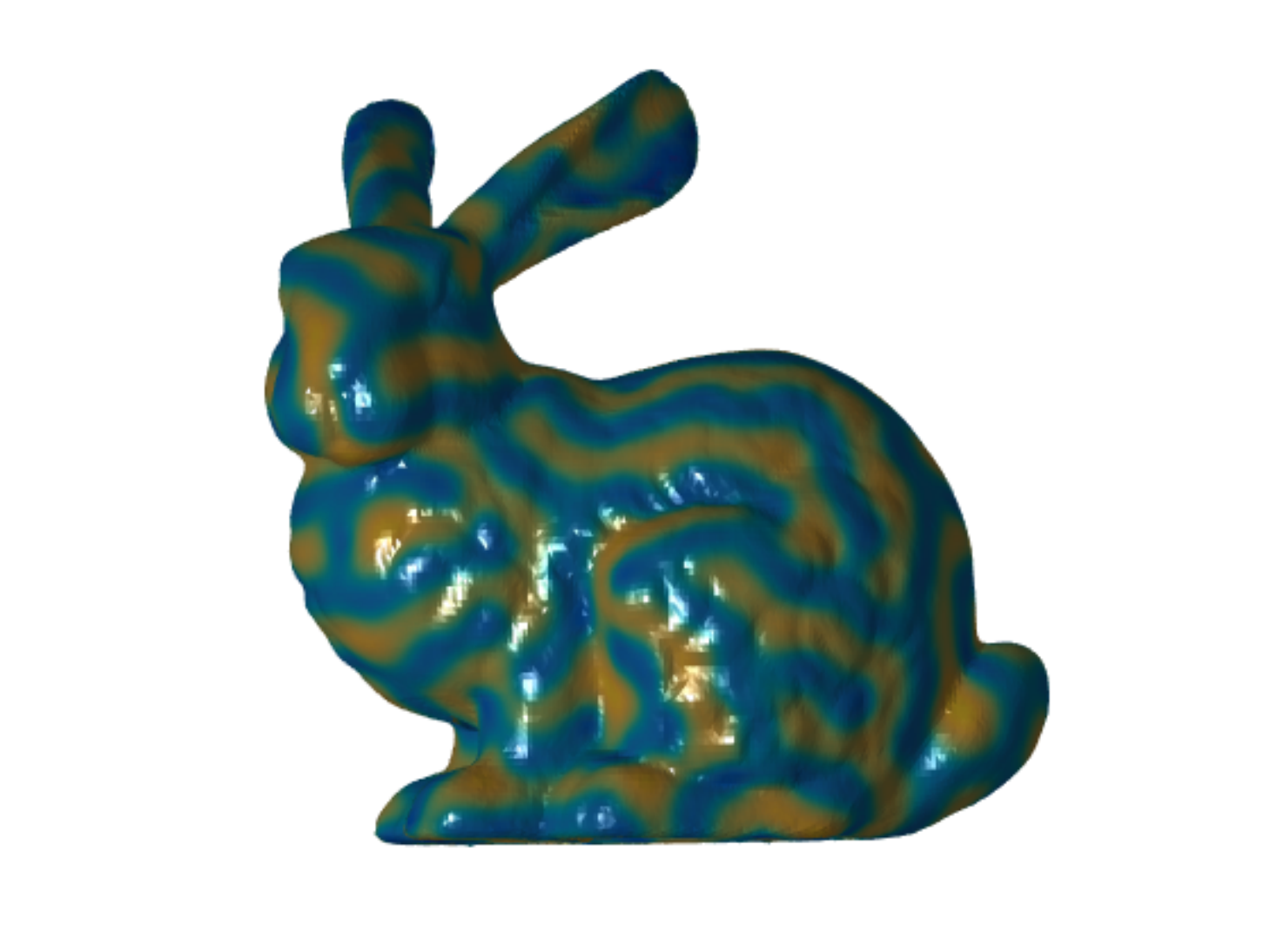}
    \caption{The solution of the Gray-Scott reaction-diffusion model for parameters $(k,F) = (0.062,0.03)$ (left) and $(k,F)=(0.06,0.037)$ (right).}
    \label{RDBunny}
\end{figure}

For parameter choices of $\mathcal{D}_u = 5\times 10^{-5}$ and $\mathcal{D}_v=2.5\times 10^{-5}$, a variety of patterns are observed as  $F$ and $k$ are varied.
Figure~\ref{RDBunny} shows two of these patterns on the surface of the Stanford Bunny \cite{turk1994zippered}. The closest point representation to the triangulated surface is calculated according to the method described in Section~\ref{section: CPCalculation}.
In this example, the final time is 15000 and the time step-size is $\Delta t = (0.1/\mathcal{D}_u) \Delta x^2$ for a spatial grid size of $\Delta x = 0.025$.
See Figure~\ref{RDBunny} for patterns arising for two different choices of the parameters $F$ and $k$ \cite{mcgough2004pattern}.

\section{Summary}\label{SummarySection}
In this paper, an explicit closest point method is introduced that uses finite differences derived from radial basis functions (RBF-FD).
In our method, an RBF-FD approximation of surface derivatives is formed using the
$m$ grid points closest to a surface point.
Localization of the computation is accomplished by computing over a tube whose radius
is obtained from the solution to the Gauss circle problem.
An advantage of our algorithm relative to the standard finite difference CPM is a reduction of the computational tube radius, leading to the reduction of the grid points in the computational domain and their corresponding closest points on the surface.
Also, higher-order schemes are easily constructed by increasing the number of points in the finite difference stencil.
When compared to RBF methods, our algorithm does not require quasi-uniform distribution of points on the surface. In addition, the repeated patterns in our computational geometry allows us to use an algorithm to invert (a small number of) collocation matrices, thereby reducing computational cost over other existing methods.
Numerical experiments are provided to validate the method for different types of PDEs on surfaces.

The RBF-FD discretization introduced in this paper solves surface PDEs using explicit time stepping methods.
Implicit RBF-FD schemes which allow for large time steps for stiff problems
are also needed, and are a focus of our current work.
Related to this, the approximation of the eigenvalues of surface operators using the RBF-CPM method
is particularly interesting (cf. \cite{macdonald2011solving}).
Another focus of our work is the
solution of PDEs on {\it moving} surfaces.
In moving closest point representations, grid node deactivation may occur  \cite{leung2009grid}.  Methods based on RBF-FD discretizations
accommodate irregular stencils and are therefore particularly attractive for such problems.
For a discussion on the issue of grid node deactivation, and an initial method
using the original closest point method, see \cite{petras2016pdes}.
\section*{Acknowledgements}
The first and third authors were partially supported by an NSERC Canada grant (RGPIN 227823).
The first author was partially supported by the Basque Government through the BERC 2014-2017 program and by Spanish Ministry of Economy and Competitiveness MINECO through BCAM Severo Ochoa excellence accreditation SEV-2013-0323 and through project MTM2015-69992-R  BELEMET. The second author was partially supported by a Hong Kong Research Grant Council GRF Grant (HKBU 11528205) and a Hong Kong Baptist University FRG Grant.

\bibliographystyle{elsarticle-num-names}
\bibliography{draftBibliography3}

\end{document}